\documentclass[10pt]{article}

\usepackage{amssymb,latexsym}
\usepackage{epsfig}
\usepackage{eufrak}
\usepackage{amsmath}
\usepackage{mathrsfs}
\usepackage{color}

\newtheorem{theorem}{Theorem}
\newtheorem{lemma}{Lemma}
\newtheorem{corollary}{Corollary}

\newcommand{\be}{\begin{equation}}
\newcommand{\ee}{\end{equation}}
\newcommand{\bee}{\begin{eqnarray*}}
\newcommand{\eee}{\end{eqnarray*}}
\newcommand{\bel}{\begin{eqnarray}}
\newcommand{\eel}{\end{eqnarray}}
\newcommand{\bec}{\begin{cases}}
\newcommand{\eec}{\end{cases}}
\newcommand{\bem}{\begin{bmatrix}}
\newcommand{\eem}{\end{bmatrix}}

\newcommand{\la}{\label}
\newcommand{\li}{\left}
\newcommand{\ri}{\right}

\newcommand{\ovl}{\overline}

\newcommand{\vep}{\varepsilon}
\newcommand{\lm}{\lambda}

\newcommand{\si}{\sigma}

\newcommand{\ga}{\gamma}

\newcommand{\vse}{\vartheta}
\newcommand{\se}{\theta}

\newcommand{\ze}{\zeta}
\newcommand{\al}{\alpha}

\newcommand{\Om}{\Omega}

\newcommand{\f}{\frac}
\newcommand{\sq}{\sqrt}
\newcommand{\cd}{\cdots}

\newcommand{\qu}{\quad}
\newcommand{\qqu}{\qquad}
\newcommand{\fa}{\forall}

\newcommand{\mscr}{\mathscr}
\newcommand{\mcal}{\mathcal}

\newcommand{\bb}{\mathbb}

\newcommand{\mrm}{\mathrm}
\newcommand{\bs}{\boldsymbol}

\newcommand{\tx}{\text}

\newcommand{\iy}{\infty}

\newcommand{\bed}{\begin{description}}
\newcommand{\eed}{\end{description}}
\newcommand{\bei}{\begin{itemize}}
\newcommand{\eei}{\end{itemize}}
\newcommand{\ben}{\begin{enumerate}}
\newcommand{\een}{\end{enumerate}}
\newcommand{\bib}{\bibitem}
\newcommand{\beL}{\begin{lemma}}
\newcommand{\eeL}{\end{lemma}}
\newcommand{\beT}{\begin{theorem}}
\newcommand{\eeT}{\end{theorem}}
\newcommand{\beC}{\begin{corollary}}
\newcommand{\eeC}{\end{corollary}}

\newcommand{\bpf}{\begin{pf}}
\newcommand{\epf}{\end{pf}}
\newcommand{\bsk}{\bigskip}

\newcommand{\pfbox}{\hfill\mbox{$\Box$}}

\newenvironment{pf}{\paragraph*{Proof{\rm.}}}{\pfbox\bigskip}

\begin{document}

%\title{{\bf Multidimensional Probabilistic Inequalities with Applications to Control of Uncertain Systems} }

\title{{\bf Uncertainty Inference with Applications to Control and Decision} }

\author{Xinjia Chen\\
Department of Engineering  Technology\\
 Northwestern State
University, Natchitoches, LA 71497\\
Email: chenx@nsula.edu \quad Tel: (318)357-5521  \quad Fax:
(318)357-6145}

\date{ }

\maketitle

\begin{abstract}

In many areas of engineering and sciences, decision rules and
control strategies are usually designed based on nominal values of
relevant system parameters. To ensure that a control strategy or
decision rule will work properly when the relevant parameters vary
within certain range, it is crucial to investigate how the
performance measure is affected by the variation of system
parameters. In this paper, we demonstrate that such issue boils down
to the study of the variation of functions of uncertainty. Motivated
by this vision, we propose a general theory for inferring function
of uncertainties. By virtue of such theory, we investigate
concentration phenomenon of random vectors. We derive uniform
exponential inequalities and multidimensional probabilistic
inequalities for random vectors, which are substantially tighter as
compared to existing ones. The probabilistic inequalities are
applied to investigate the performance of control systems with real
parametric uncertainty. It is demonstrated much more useful insights
of control systems can be obtained. Moreover, the probabilistic
inequalities offer performance analysis in a significantly less
conservative way as compared to the classical deterministic
worst-case method.

\end{abstract}

%\tableofcontents

\section{Introduction}

Decision and control are frequent problems of many areas of
engineering and sciences. In general, the object that we are facing
and need to design proper control strategy or decision rule can be
viewed as a system.  In most cases, we don't have complete
information about such system. In order to avoid system failure, it
is an essential task to evaluate the performance of the systems
affected by uncertainty \cite{Dorf, Franklin}. Existing methods for
performance evaluation of uncertain systems are based on two
completely different paradigms. The first paradigm is to treat
uncertainty as deterministic bounded parameters \cite{Barmishbook,
Zhou}. The performance analysis is to seek the worst-case scenario.
This approach can be unduely conservative. Moreover, the resultant
computational complexity can be NP hard. The second paradigm is to
evaluate system performance by assuming some typical distribution
for the underlying uncertainty \cite{Astrom, DTCRC}. This approach
can be conducted with Monte Carlo simulation. The computational
complexity can be shown to be independent of the problem size. The
major issue of such paradigm is that the assumed distribution may be
significantly different from the actual distribution of the
underlying uncertainty.  Consequently, the resultant insight from
the Monte Carlo simulation can be fairy misleading.

Actually, in the analysis and design of control strategies and
decision rules, due to experimental or cognitive limitations, we
only have limited information about the uncertainty affecting the
systems \cite{ChenGu, Ljung}. Motivated by this situation, we
advocate to analyze system performance based on the limited
available information. Specifically, we represent such information
by constraints of the mathematical expectation of functions of
uncertainty. The performance measure of systems is expressed as the
mathematical expectation of certain functions of uncertainty.
Consequently, the range of such expected value is a good indicator
of the performance of the associated system.  In this way, we
establish a close connection between probabilistic inequalities and
the analysis and design of control and decision. More formally, the
general problem can formulated as follows. Let $X$ be a random
vector representing uncertainty affecting the systems. Let $\bs{f}
(.)$ be a function of the uncertainty and $\mcal{D}$ be a domain in
the Euclidean space such that $\bb{E} [ \bs{f} (X) ] \in \mcal{D}$.
Let $g (X)$ denote the performance of the system. It is desirable to
determine the range of $\bb{E} [ g (X)]$.  This formulation
accommodate a wide range of problems on performance analysis of
control systems as special cases.  A familiar problem is the robust
stability of uncertain system. Within this general framework, we
derive tight bounds for $\bb{E} [ g (X)]$, which can be evaluated by
computational techniques such as linear programming embedded with
gradient search \cite{Bazaraa} and global optimization techniques
such as branch and bound algorithm \cite{Land}.

The remainder of the paper is organized as follows.  In Section 2,
we propose a general approach for inferring uncertainty.  Such
approach is based on a probabilistic characterization of convex
sets. In Section 3, we apply the proposed theory of inferring
uncertainty to investigate concentration phenomena frequently
encountered in uncertain systems.  We shall first establish uniform
exponential concentration inequalities.  Afterward, multidimensional
probabilistic inequalities are developed which are useful for
analysis of control systems.  In Section 4, we apply the
probabilistic theory to analyze the stability of control systems
affected by parametric uncertainty. Section 5 is the conclusion.
Most proofs are given in Appendices.

In this paper, we shall use the following notations.  The set of
real numbers is denoted by $\bb{R}$.  The set of nonnegative real
numbers is denoted by $\bb{R}^+$. The $d$-dimensional Euclidean
space is denoted by $\bb{R}^d$.   The set of positive integers is
denoted by $\bb{N}$. The Euclidean norm is denoted by $||. ||$. The
diameter of $S \subseteq \bb{R}^d$ is defined as $\sup \{ || x - y
||: \; x \in S, \; y \in S \}$.  The supremum of an empty set is
defined as $0$. The set minus operation is denoted by $\setminus$.
Let $(\Om, \mscr{F}, \Pr )$ denote the probability space.  The
mathematical expectation of random vector $X$ is denoted by $\bb{E}
[ X ]$.  A zero-mean random vector is a random vector such that all
the elements of its expected value are zero.

Let $X$ be a discrete random vector in $\bb{R}^d$. A vector $x$ in
$\bb{R}^d$ is said to be a possible value of the discrete random
vector if $\Pr \{ X = x \} > 0$.   That is, a vector $x$ is said to
be a possible value of a discrete random vector if $x$ is assumed by
the discrete random vector with a positive probability.

The support of a random variable $X$ in $\bb{R}^d$ is defined as the
set whose complement consists of points in $\bb{R}^d$ with zero
probability density.  We use the abbreviation ``i.i.d.'' for
``independent and identically distributed''. The first and second
derivatives of function $\psi (s)$ is denoted by $\psi^\prime (s)$
and $\psi^{\prime \prime} (s)$, respectively.  We use the big O
notation $f(x) = O (g (x))$ as $x \to a$ in the sense that
$\limsup_{x \to a} \li | \f{ f (x) }{ g (x)} \ri | < \iy$.  The
other notations will be made clear as we proceed.

\section{A General Theory for Inferring Uncertainty}

In this section, we shall develop a general theory for inferring
uncertainty.  To make the inference more realistic, we avoid the
assumption that the exact distribution of uncertainty is known.  We
shall demonstrate that a unified theory of inference can be
established upon  a stochastic characteristic of convex sets.

\subsection{ A Stochastic Characteristic of Convex Sets}

Our investigation indicates that if a set in a finite-dimensional
Euclidean space is convex, then the set contains the expectation of
any random vector almost surely contained by the set.  More
formally, we have established the following result.

\beT

\la{fundamental}

If $\mscr{D}$ is a convex set in $\bb{R}^n$, then $\bb{E} [
\bs{\mcal{X}} ] \in \mscr{D}$ holds for any random vector
$\bs{\mcal{X}}$ such that $\Pr \{ \bs{\mcal{X}} \in \mscr{D} \} = 1$
and that $\bb{E} [ \bs{\mcal{X}} ] $ exists.
  \eeT

Theorem \ref{fundamental} is established in \cite{Chen2015}.   The
converse of Theorem \ref{fundamental} asserts that if $\mscr{D}$ is
a set in $\bb{R}^n$ such that $\bb{E} [ \bs{\mcal{X}} ] \in
\mscr{D}$ holds for any random vector $\bs{\mcal{X}}$ such that $\Pr
\{ \bs{\mcal{X}}  \in \mscr{D} \} = 1$ and that $\bb{E} [
\bs{\mcal{X}} ] $ exists, then $\mscr{D}$ is convex. This assertion
is well known and is a direct consequence of the definition of a
convex set.

Theorem \ref{fundamental} immediately implies Jensen's inequality.
To see this, note that if a function is convex, then its epigraph,
the region above its graph, is a convex set. Hence, if $f$ is a
convex function, then for any random variable $X$, since $(X, \;
f(X))$ is contained by the epigraph of $f$, it follows from Theorem
\ref{fundamental} that $(\bb{E}[ X ], \; \bb{E} [ f(X) ] )$ is
contained by its epigraph. This implies that $\bb{E} [ f(X) ] \geq f
(\bb{E} [ X ] )$ by the notion of epigraph.

The following result is due to Isii \cite{Isii}.

\beT

\la{theoremIsii}

Let $\mscr{X}$ be a family of random vectors in $\bb{R}^d$ such that
\[
\Pr \{ X \in \mscr{A} \} = 1, \qqu   \bb{E} [ \bs{f} (X) ] =
\bs{\mu} \in \bb{R}^k  \qu \tx{for each $X \in \mscr{X}$},
\]
where $\mscr{A}$ is a subset of $\bb{R}^d$ and $\bs{f} (x)$ is a
function assuming values in $\bb{R}^k$ for $x \in \mscr{A}$.  Let
$g(x)$ be real-valued function of $x \in \mscr{A}$  such that
$\bb{E} [ g (X) ]$ exists for each $X \in \mscr{X}$.  Then,
\[
\sup_{ X \in \mscr{X} } \bb{E} [ g (X) ] =  \sup_{ Y \in \mscr{Y} }
\bb{E} [ g (Y) ],
\]
where \[ \mscr{Y} = \{ Y \in \mscr{X}: \tx{ Y is a discrete random
vector with at most $k +1$ distinct possible values} \}.
\]
\eeT

This result is correct.  However,  in his original proof, Isii made
a mistake by using an incorrect probability measure in mathematical
induction (see, \cite[Lemma 2, page 191--192]{Isii}).

\bsk

In many applications, because of incomplete information, the
equality $\bb{E} [ \bs{f} (X) ] = \bs{\mu}$ is hard to satisfy. For
example, in many cases, we may not know the exact value of the
moment of a random variable.  We only have its range.   Hence, to
infer uncertainty in the most general setting, we propose to
represent the incomplete information by the constraint
\[
\bb{E} [ \bs{f} (X) ] \in \mscr{B}, \]
 where  $\mscr{B}$ is a subset of $\bb{R}^k$.  In this framework, we have the following result.

\beT

\la{theorem88388chen}

Let $X$ be a random vector in $\bb{R}^d$ such that $\Pr \{ X \in
\mscr{A} \} = 1$ and $\bb{E} [ \bs{f} (X) ] \in \mscr{B}$,  where
$\mscr{A}$ is a subset of $\bb{R}^d$,  $\mscr{B}$ is a subset of
$\bb{R}^k$, and $\bs{f} (x)$ is a function assuming values in
$\bb{R}^k$ for $x \in \mscr{A}$.    Let $g(x)$ be a real-valued
function of $x \in \mscr{A}$ such that $\bb{E} [ g (X) ]$ exists.
Then,
\[
\bb{E} [ g (X) ] \leq \sup_{ Y \in \mscr{Y} } \bb{E} [ g (Y) ],
\]
where $\mscr{Y}$ is the family of discrete random vectors in
$\bb{R}^d$ such that for each $Y \in \mscr{Y}$,
\[
\Pr \{ Y \in \mscr{A} \} = 1, \qqu   \bb{E} [ \bs{f} (Y) ] \in
\mscr{B},
\]
and $Y$ has at most $k +1$ distinct possible values. \eeT

See Appendix \ref{theorem88388chenapp} for a proof.   Making use of
Theorem \ref{theorem88388chen}, we have the following result.

\beT

\la{theorem883}

Let $\mscr{X}$ be a family of random vectors in $\bb{R}^d$ such that
\[
\Pr \{ X \in \mscr{A} \} = 1, \qqu   \bb{E} [ \bs{f} (X) ] \in
\mscr{B} \qu \tx{for each $X \in \mscr{X}$},
\]
where $\mscr{A}$ is a subset of $\bb{R}^d$, $\mscr{B}$ is a subset
of $\bb{R}^k$, and $\bs{f} (x)$ is a function assuming values in
$\bb{R}^k$ for $x \in \mscr{A}$.  Let $g(x)$ be real-valued function
of $x \in \mscr{A}$  such that $\bb{E} [ g (X) ]$ exists for each $X
\in \mscr{X}$.  Then,
\[
\sup_{ X \in \mscr{X} } \bb{E} [ g (X) ] =  \sup_{ Y \in \mscr{Y} }
\bb{E} [ g (Y) ],
\]
where \[ \mscr{Y} = \{ Y \in \mscr{X}: \tx{ Y is a discrete random
vector with at most $k +1$ distinct possible values} \}.
\]
\eeT

Theorem \ref{theorem883} can be shown as follows.

By the assumption that $\bb{E} [ g (X) ]$ exists for each $X \in
\mscr{X}$, according to Theorem \ref{theorem88388chen}, we have that
\[
\bb{E} [ g (X) ] \leq  \sup_{ Y \in \mscr{Y} } \bb{E} [ g (Y) ]
\]
for each $X \in \mscr{X}$.  Thus,
\[
\sup_{ X \in \mscr{X} } \bb{E} [ g (X) ] \leq  \sup_{ Y \in \mscr{Y}
} \bb{E} [ g (Y) ],
\]
On the other hand, since $\mscr{Y}$ is a subset of $\mscr{X}$, it
must be true that
\[
\sup_{ X \in \mscr{X} } \bb{E} [ g (X) ] \geq  \sup_{ Y \in \mscr{Y}
} \bb{E} [ g (Y) ].
\]
So, the theorem must be true.

\bsk

According to Theorem \ref{theorem883}, we have {\small \bee \sup_{ Y
\in \mscr{Y} } \bb{E} [ g (Y) ]  =  \sup \li \{ \sum_{\ell = 1}^{k +
1} \se_\ell g (y_\ell) : \; \se_\ell \geq 0 \; \tx{and} \; y_\ell
\in \mscr{A} \; \tx{for} \; \ell = 1, \cd, k + 1, \qu
\sum_{\ell=1}^{k + 1} \se_\ell = 1, \qu \sum_{\ell = 1}^{k + 1}
\se_\ell \bs{f} ( y_\ell ) \in \mscr{B} \ri \}, \eee} which can be
computed by linear programming embedded with gradient search
\cite{Bazaraa}, and branch and bound method \cite{Land}.

For the important case that $g(.)$ is an indicator function, we have
the following result.

\beT

\la{VIPProb}

Let $\mscr{X}$ be a family of random vectors in $\bb{R}^d$ such that
$\Pr \{ X \in \mscr{A} \} = 1$ and $\bb{E} [ \bs{f} (X) ] \in
\mscr{B}$ for each $X \in \mscr{X}$,  where $\mscr{A}$ is a subset
of $\bb{R}^d$,   $\mscr{B}$ is a subset of $\bb{R}^k$, and $\bs{f}
(x)$ is a function assuming values in $\bb{R}^k$ for $x \in
\mscr{A}$. Then, $\sup_{ X \in \mscr{X} }  \Pr \{ X \in \mscr{C}  \}
= \max \{ P_i: 1 \leq i \leq k + 1  \}$ for any subset $\mscr{C}$ of
$\mscr{A}$, where \bee P_i  & = & \sup \li \{  \sum_{\ell = 1}^i
\se_\ell: \; \se_\ell \geq 0  \; \tx{for} \; 1 \leq \ell \leq  k +
1, \qu y_\ell \in \mscr{C} \; \tx{for} \; 1 \leq \ell \leq i,
 \qu  y_\ell \in \mscr{A} \setminus \mscr{C} \; \tx{for} \; i < \ell \leq k + 1, \ri. \\
&  &  \qu \qu \qu \qu \qu \qu \li.  \sum_{\ell = 1 }^{k + 1}
\se_\ell = 1, \qu \sum_{\ell = 1}^{k + 1} \se_\ell \bs{f} ( y_\ell )
\in \mscr{B} \ri \} \eee for $i = 1, \cd, k + 1$. \eeT

See Appendix \ref{VIPProbapp} for a proof.

\bsk

Theorem \ref{VIPProb} can be applied to compute bounds for the
probability that a systems fails to satisfy pre-specified
requirements based on limited information of uncertainty. The bounds
can be obtained by linear programming embedded with gradient search,
and the branch and bound method. A demonstration of the application
of this theorem is given in Section 4.

\subsection{Minimum-Range Random Variable Under Moment Constraints}

Making use of Theorem \ref{theorem883}, we have the following
result.

\beT

Let $Z$ be a zero-mean random variable in $\bb{R}$  such that \be
\la{constraint}  \bb{E} [ Z^k ] \geq 1 \qu \tx{for} \qu k \geq 2.
\ee Define $L_Z = \sup \{ u \in \bb{R}: \Pr \{ Z \geq u \} = 1 \}$
and $U_Z = \inf \{ v \in \bb{R}: \Pr \{ Z \leq v \} = 1 \}$.
 Then, $U_Z - L_Z \geq \sq{5}$.  In particular, (\ref{constraint}) holds and $U_Z - L_Z = \sq{5}$ if
$Z$ is a random variable such that $\Pr \li \{ Z = \varphi \ri \} =
\f{1}{\sq{5} \; \varphi}$ and $\Pr \li \{ Z = - \f{1}{\varphi} \ri
\} = \f{\varphi}{\sq{5}}$, where $\varphi = \f{ 1 + \sq{5} }{2}$ is
the golden ratio.

\eeT

Making use of Theorem \ref{theorem883}, we have the following
result.

 \beT

Let $Z$ be a zero-mean random variable in $\bb{R}$  such that \be
\la{constraintb} \bb{E} [ Z^2 ] = 1, \qqu \bb{E} [ Z^k ] \geq 1 \qu
\tx{for} \qu k \geq 3. \ee Define $L_Z = \sup \{ u \in \bb{R}: \Pr
\{ Z \geq u \} = 1 \}$ and $U_Z = \inf \{ v \in \bb{R}: \Pr \{ Z
\leq v \} = 1 \}$.
 Then, $\max (  U_Z , \; | L_Z | ) \geq \varphi$, where $\varphi = \f{ 1 + \sq{5} }{2}$ is the golden ratio.
In particular, (\ref{constraintb}) holds and $\max (  U_Z , \; | L_Z
| ) = \varphi$ if $Z$ is a random variable such that $\Pr \li \{ Z =
\varphi \ri \} = \f{1}{\sq{5} \; \varphi}$ and $\Pr \li \{ Z = -
\f{1}{\varphi} \ri \} = \f{\varphi}{\sq{5}}$.

\eeT

\section{Concentration Phenomena in Euclidean Space}

In many applications, uncertainties can be represented as random
vectors in Euclidean space.  Consequently, useful insight of the
impact of uncertainty to control and decision may be obtained  by
investigating the concentration phenomena of the relevant random
vectors. In the sequel, we shall develop concentration inequalities
for random vectors, which generalize Chernoff-Hoeffding inequalities
\cite{Chernoff, Hoeffding}.  For that purpose, we shall first
propose a unified approach for deriving exponential inequalities
which uniformly hold for all values of time for stochastic
processes.

\subsection{Uniform Exponential Inequalities}

The following results provide a unified method for deriving uniform
exponential inequalities for real-valued stochastic processes.

\beT  \la{Them9} [Chen (2012)] Let $\mcal{V}_t$ be a non-negative,
right-continuous function of $t \in [0, \iy )$. Let $\{ X_t, \; t
\in \bb{R}^+ \}$ be a right-continuous stochastic process such that
$\bb{E} [ \exp ( s (X_{t^\prime} - X_t)  ) \mid \mscr{F}_t ] \leq
\exp ( (\mcal{V}_{t^\prime} - \mcal{V}_t) \varphi(s) )$ almost
surely for arbitrary $t^\prime \geq t \geq 0$ and $s \in (0, b)$,
where $b$ is a positive number or infinity,  $\varphi(s)$ is a
non-negative function of $s \in (0, b)$, and $\mscr{F}_t$ is the
$\si$-algebra generated by $\{ X_{t^\prime}, \; 0 \leq t^\prime \leq
t \}$. Let $\tau > 0$ and $\ga > 0$. Then, {\small \be \Pr \li \{
\sup_{t > 0 } \li [ X_t - X_0 - \ga \mcal{V}_\tau  - \f{ \varphi(s)
}{s} (\mcal{V}_t - \mcal{V}_\tau) \ri ] \geq  0 \ri \} \leq \li [
\exp \li ( \varphi(s) - \ga s \ri ) \ri ]^{\mcal{V}_\tau} \qqu \fa s
\in (0, b). \la{genb12} \ee}  In particular, if $\{ s \in (0, b):
\varphi(s) \leq \ga s \}$ is nonempty and the infimum of $\varphi(s)
- \ga s$ with respect to $s \in (0, b)$ is attained at $\zeta \in
(0, b)$, then {\small \be \Pr \li \{ \sup_{t > 0} \li [  X_t - X_0 -
\ga \mcal{V}_\tau - \f{ \varphi(\zeta) }{\zeta} (\mcal{V}_t -
\mcal{V}_\tau) \ri ] \geq 0 \ri \} \leq [ \exp ( \varphi(\zeta) -
\ga \zeta ) ]^{\mcal{V}_\tau} \leq 1, \la{discoverb12} \ee} and $0
\leq \f{ \varphi(\zeta) }{\zeta} \leq \ga$. \eeT

Theorem \ref{Them9} is established in \cite{ChenSPIE2012, Chen2012}.
A proof is reproduced in  Appendix \ref{Them9_app}.   More
generally, we have the following results.

\beT  \la{Them988}  Let $\{ \mcal{V} (s, t), \; t \in \bb{R}^+ \}$
be a real-valued stochastic process parameterized by $s \in (0, b)$,
where $b$ is a positive number or infinity.  Let  $\{ X_t, \; t \in
\bb{R}^+ \}$ be a real-valued stochastic process with $X_0 = 0$. Let
$\{ \mcal{Z} (s, t), \; t \in \bb{R}^+ \}$ be a right-continuous
 supermartingale,  which is parameterized by $s \in (0, b)$ and adapted to
 the natural filtration generated by $\{ \mcal{V} (s, t), \; t \in
\bb{R}^+ \}$ and $\{ X_t, \; t \in \bb{R}^+ \}$ such that for all $s
\in (0, b)$,
\[
\bb{E} [ \mcal{Z} (s, 0) ] \leq 1 \qu \tx{and} \qu \exp ( s X_t  -
\mcal{V} (s, t) ) \leq \mcal{Z} (s, t) \qu \tx{almost surely for all
$t \in \bb{R}^+$}.
\]
Let $\ga$ be a real number and $g (s)$ be a function of $s \in (0,
b)$. Then, \be \la{gen88338} \Pr \li \{ \sup_{t > 0} \li [ X_t - \ga
- \f{ \mcal{V} (s, t) - g (s) } {s} \ri ] \geq 0  \ri \} \leq \exp (
g (s) - \ga s ) \qu \tx{for all $s \in (0, b)$.} \ee In particular,
the following assertions hold:

\bed

\item (I)  If the infimum of $g (s ) - \ga s$ with respect to
$s \in (0, b)$ is attained at $\zeta \in (0, b)$, then {\small \[
\Pr \li \{ \sup_{t > 0} \li [ X_t - \ga - \f{ \mcal{V} (\zeta, t) -
g (\zeta) } {\zeta} \ri ] \geq 0  \ri \} \leq \exp ( g (\zeta) - \ga
\zeta ).
\]}

\item (II) If $\mcal{V}(s, t)$ is a deterministic function of $s \in
(0, b)$ and $t \in \bb{R}^+$, then {\small \[ \Pr \li \{ \sup_{t >
0} \li [ X_t - \ga - \f{ \mcal{V} (s, t) - \mcal{V}(s, \tau ) } {s}
\ri ] \geq 0  \ri \} \leq  \exp ( \mcal{V}(s, \tau ) - \ga s ) \] }
for all $s \in (0, b)$ and $\tau \in \bb{R}^+$.

\item (III) If $\mcal{V}(s, t) = \varphi(s) V_t + \ln C$, where $C$ is a positive constant,  $\varphi(s)$
is a deterministic function of $s \in (0, b)$,  and $\{ V_t, \; t
\in \bb{R}^+ \}$ is a deterministic or stochastic process, then
{\small
\[ \Pr \li \{ \sup_{t > 0} \li [ X_t - \ga - \f{ \varphi (s) } {s}
(V_t - m) \ri ] \geq 0 \ri \} \leq C \exp ( m \varphi (s)  - \ga s )
\] } for all $s \in (0, b)$ and $m \in \bb{R}$.
\eed

\eeT

Theorem \ref{Them988} is established in \cite{ChenSPIE2020} and
presented in SPIE Conference in April 2020.  It should be noted that
if $\varphi (s)$ has the characteristic of a cumulant-generating
function, then the assertion (III) of Theorem \ref{Them988} can be
applied to deduce Theorem 1(b) of \cite{Howard}.

\bsk

To prove Theorem \ref{Them988}, note that for all $s \in (0, b)$,
{\small \bee \Pr \li \{ \sup_{t > 0} \li [ X_t - \ga - \f{ \mcal{V}
(s, t) - g (s) } {s} \ri ] \geq 0  \ri \}  & = & \Pr \li \{ \sup_{t
> 0} s \li [ X_t - \ga - \f{ \mcal{V} (s, t) - g (s) } {s}
\ri ] \geq 0  \ri \} \\
& = & \Pr  \li \{ \sup_{t > 0} \li [ s X_t - \mcal{V} (s, t) \ri ]
\geq \ga s - g (s) \ri  \} \\
& = & \Pr \li \{ \sup_{t > 0} \exp \li ( s X_t - \mcal{V} (s, t) \ri
) \geq \exp (\ga s - g (s) ) \ri \}  \\
& \leq  & \Pr \li \{ \sup_{t
> 0} \mcal{Z} (s, t)  \geq \exp ( \ga s - g (s) ) \ri \}. \eee } By the
supermartingale inequality, we have {\small \[ \Pr \li \{ \sup_{t >
0} \li [ X_t - \ga - \f{ \mcal{V} (s, t) - g (s) } {s} \ri ] \geq 0
\ri \} \leq \f{ \bb{E} [ \mcal{Z} (s, 0) ] }{ \exp ( \ga s - g (s) )
} \leq \f{ 1 }{ \exp ( \ga s - g (s) ) } = \exp (g (s) - \ga s )
\] }
for all $s \in (0, b)$. This proves (\ref{gen88338}), from which the
particular assertions  immediately follow.

\bsk

Theorem \ref{Them988} concerns the probability of crossing the curve
in the upward direction.  Similar results can be derived for the
probability of crossing a curve in the downward direction. Moreover,
it is possible to unify the inequalities for the probabilities of
crossing curves in both upward and downward  directions  by the
following results.

\beT  \la{Them98}  Let $\{ \mcal{V} (s, t), \; t \in \bb{R}^+ \}$ be
a real-valued stochastic process parameterized by $s \in \mscr{S}
\subseteq \bb{R}$.  Let $\{ X_t, \; t \in \bb{R}^+ \}$ be a
real-valued stochastic process with $X_0 = 0$. Let $\{ \mcal{Z} (s,
t), \; t \in \bb{R}^+ \}$ be a right-continuous
 supermartingale,  which is parameterized by $s \in \mscr{S}$ and adapted to
 the natural filtration generated by $\{ \mcal{V} (s, t), \; t \in
\bb{R}^+ \}$ and $\{ X_t, \; t \in \bb{R}^+ \}$ such that for all $s
\in \mscr{S}$,
\[
\bb{E} [ \mcal{Z} (s, 0) ] \leq 1 \qu \tx{and} \qu \exp ( s X_t  -
\mcal{V} (s, t) ) \leq \mcal{Z} (s, t) \qu \tx{almost surely for all
$t \in \bb{R}^+$}.
\]
Let $\ga$ be a real number and $g (s)$ be a function of $s \in
\mscr{S}$. Then, \be \la{gen8833} \Pr \li \{ \sup_{t > 0} \li [ s (
X_t - \ga ) -  \mcal{V} (s, t) + g (s) \ri ] \geq 0  \ri \} \leq
\exp ( g (s) - \ga s ) \qu \tx{for all $s \in \mscr{S}$.} \ee In
particular, the following assertions hold:

\bed

\item (I)  If the infimum of $g (s )- \ga s$ with respect to $s \in
\mscr{S}$ is attained at $\zeta \in \mscr{S}$, then {\small \[ \Pr
\li \{ \sup_{t > 0} \li [ \zeta ( X_t - \ga ) - \mcal{V} (\zeta, t)
+ g (\zeta)  \ri ] \geq 0  \ri \} \leq  \exp ( g (\zeta) - \ga \zeta
). \] }

\item (II) If $\mcal{V}(s, t)$ is a deterministic function of $s \in
\mscr{S}$ and $t \in \bb{R}^+$, then {\small \[ \Pr \li \{ \sup_{t >
0} \li [ s ( X_t - \ga ) - \mcal{V} (s, t) + \mcal{V}(s, \tau ) \ri
] \geq 0  \ri \} \leq  \exp ( \mcal{V}(s, \tau ) - \ga s ) \] } for
all $s \in \mscr{S}$ and $\tau \in \bb{R}^+$.

\item (III) If $\mcal{V}(s, t) = \varphi(s) V_t + \ln C$, where $C$ is a positive constant, $\varphi(s)$
is a deterministic function of $s \in \mscr{S}$,  and $\{ V_t, \; t
\in \bb{R}^+ \}$ is a deterministic or stochastic process, then
{\small
\[ \Pr \li \{ \sup_{t > 0} \li [ s ( X_t - \ga ) - \varphi (s) (V_t
- m) \ri ] \geq 0 \ri \} \leq C \exp ( m \varphi (s)  - \ga s ) \] }
for all $s \in \mscr{S}$ and $m \in \bb{R}$.

\eed

\eeT

See Appendix \ref{Them98app}  for a proof.

Making use of Theorem \ref{Them98}, we have the following result.

\beT \la{simple} Let $Y_1, Y_2, \cd$ be a sequence of  independent
random variables.  Define $X_n = \sum_{i =1}^n Y_i$ for $n \in
\bb{N}$. Assume that the moment generating function, $\bb{E} [ e^{s
X_n} ]$, is bounded from above by $ \exp ( \mcal{V} (s, n) )$ ,
where $\mcal{V} (s, n)$ is a function of $s \in \mscr{S} \subseteq
\bb{R}$ and $n \in \bb{N}$. Let $\ga \in \bb{R}$ and $m \in \bb{N}$.
Then, {\small
\[ \Pr \li \{ \sup_{n \in \bb{N} } \li [ s ( X_n - \ga ) - \mcal{V}
(s, n) + \mcal{V}(s, m ) \ri ] \geq 0  \ri \} \leq \exp (
\mcal{V}(s, m ) - \ga s ) \] } for all $s \in \mscr{S}$.  Specially,
if the infimum of $\mcal{V}(s, m ) - \ga s$ with respect to $s \in
\mscr{S}$ is attained at $\zeta \in \mscr{S}$, then {\small \[ \Pr
\li \{ \sup_{n \in \bb{N} } \li [ \zeta ( X_n - \ga ) - \mcal{V}
(\zeta, n) + \mcal{V}(\zeta, m ) \ri ] \geq 0  \ri \} \leq \exp (
\mcal{V}(\zeta, m ) - \ga \zeta ).  \]}

\eeT

See Appendix \ref{simpleapp} for a proof.

As an immediate application of Theorem \ref{simple}, we have the
following result.

\beT \la{special} Let $Y, Y_1, Y_2, \cd$ be a sequence of i.i.d.
 random variables.  Define $X_n = \sum_{i =1}^n Y_i$ for
$n \in \bb{N}$. Assume that the moment generating function, $\bb{E}
[ e^{s Y} ]$, is bounded from above by $ \exp ( \varphi (s) )$ ,
where $\varphi (s)$ is a function of $s \in \mscr{S} \subseteq
\bb{R}$. Let $\se \in \bb{R}$ and $m \in \bb{N}$. Then, {\small
\[ \Pr \li \{ \sup_{n \in \bb{N} } \li [ s ( X_n - m \se) - \varphi
(s) (n - m )  \ri ] \geq 0  \ri \} \leq \li [ \exp ( \varphi (s)  -
\se s ) \ri ]^m
\] } for all $s \in \mscr{S}$. Specially, if the infimum of
$\varphi (s)  - \se s $ with respect to $s \in \mscr{S}$ is attained
at $\zeta \in \mscr{S}$, then  {\small
\[ \Pr \li \{ \sup_{n \in \bb{N} } \li [ \zeta ( X_n - m \se) - \varphi
(\zeta) (n - m )  \ri ] \geq 0  \ri \} \leq \li [ \exp ( \varphi
(\zeta) - \se \zeta ) \ri ]^m.
\] }

\eeT

It is interesting to investigate the asymptotic structure of the
uniform exponential inequality as the magnitude of deviation tends
to $0$.  For this purpose, we have the following results.

\beT  \la{Aspstruecure88689} Let $\mcal{V}_t$ be a non-negative,
right-continuous function of $t \in [0, \iy )$. Let $\{ X_t, \; t
\in \bb{R}^+ \}$ be a right-continuous stochastic process with $X_0
= 0$ such that $\{ \exp ( s X_t - \varphi (s) \mcal{V}_t  ), \; t
\in \bb{R}^+ \} $ is a supermartingale for $s \in (a, b)$  with $a <
0 < b$, where $\varphi(s)$ is a convex function of $s \in (a, b)$
such that $\varphi(0) = \varphi^\prime (0) = 0, \; \varphi^{\prime
\prime} (0) = \al > 0$,   and has third derivative which is
continuous at a neighborhood of $0$. Let $\tau
> 0$. Then, for any real number $\vep$, the inequality {\small \be
\Pr \li \{ \sup_{t > 0 } \li [ s( X_t- \vep \mcal{V}_\tau ) -
\varphi(s) (\mcal{V}_t - \mcal{V}_\tau) \ri ] \geq 0 \ri \} \leq \li
[ \exp \li ( \varphi(s) - \vep s \ri ) \ri ]^{\mcal{V}_\tau}
\la{firstinequality689} \ee} holds for all $s \in (a, b)$. In
particular,  for $\vep$ satisfying $\lim_{s \downarrow a} \f{
\varphi (s) }{s} < \vep < \lim_{s \uparrow b} \f{ \varphi (s) }{s}$,
{\small \be \Pr \li \{ \sup_{t > 0 } \li [ \ze ( X_t- \vep
\mcal{V}_\tau ) - \varphi(\ze) (\mcal{V}_t - \mcal{V}_\tau) \ri ]
\geq 0 \ri \} \leq \li [ \exp \li ( \varphi(\ze) - \vep \ze \ri )
\ri ]^{\mcal{V}_\tau} = \inf_{s \in (a, b)} \li [ \exp \li (
\varphi(s) - \vep s \ri ) \ri ]^{\mcal{V}_\tau},
\la{secondinequality689} \ee} where $\ze$ is the unique root of the
equation $\varphi^\prime (s) = \vep$ with respect to $s \in (a, b)$,
\[
\li [   \exp  ( \varphi (\zeta) - \vep \zeta ) \ri ]^{\mcal{V}_\tau}
 =  \li [  1 +  O (\vep^3)  \ri ] \li [  \exp \li (  -
\f{\vep^2}{2 \al} \ri )  \ri ]^{\mcal{V}_\tau}  \qqu \tx{and} \qqu
\f{ \varphi (\zeta) }{\zeta} =   \f{\vep}{2}  \li [ 1 + O (\vep) \ri
]
\]
as $\vep \to 0$.  \eeT

See Appendix \ref{Aspstruecure88689app} for a proof.

For i.i.d. random variables, we have the following results.

 \beT  \la{Aspstruecure}

 Let $Y, Y_1, Y_2, \cd$ be a sequence of i.i.d. random variables such that $a = \inf \{ s < 0: \bb{E} [ e^{s Y} ] < \iy \}
 < 0$ and  $b = \sup \{ s > 0: \bb{E} [ e^{s Y} ] < \iy \} >
 0$.  Let $\varphi (s) = \ln \bb{E} [ e^{s
 (Y - \mu)} ]$ for $s \in (a, b)$.  Let $\mu = \bb{E} [ Y], \; \si^2 = \bb{E} [ | Y - \mu |^2
 ]$,
 and $\nu = \bb{E} [ (Y - \mu)^3]$.  Define $X_n =\sum_{i = 1}^n (  Y_i - \mu )$ for $n \in \bb{N}$.  Assume that $\si > 0$.
 Then, for any positive integer
 $m$ and real number $\vep$ such that $\lim_{ s \downarrow a }  \f{ \varphi
 (s)} {s} < \vep < \lim_{ s \uparrow b }  \f{ \varphi
 (s)} {s}$,  the inequality
 \be
 \la{firstinequality}
\Pr \li \{  \sup_{n \in \bb{N} }  \li [  s ( X_n  - m \vep ) -
\varphi (s) (n - m)  \ri ]  \geq 0\ri \} \leq \li [   \exp ( \varphi
(s) - \vep s )  \ri ]^m
 \ee
holds for all $s \in (a, b)$.  In particular,  \be
\la{secondinequality}
 \Pr \li \{ \sup_{n \in \bb{N} }  \li [  \zeta (
X_n - m \vep ) - \varphi (\zeta) (n - m) \ri ] \geq 0\ri \} \leq \li
[ \exp  ( \varphi (\zeta) - \vep \zeta ) \ri ]^m,
 \ee
 where $\zeta$ is the unique root of the equation $\varphi^\prime (s)
 = \vep$ with respect to $s \in (a, b)$,
\[
\li [   \exp  ( \varphi (\zeta) - \vep \zeta ) \ri ]^m  = \li [ 1 +
O (\vep^4) \ri ] \li [ \exp \li ( -
\f{\vep^2}{2 \si^2} + \f{\nu \vep^3}{6 \si^6}  \ri )  \ri ]^m\\
 =  \li [  1 +  O (\vep^3)  \ri ] \li [  \exp \li (  -
\f{\vep^2}{2 \si^2} \ri )  \ri ]^m
\]
and \[ \f{ \varphi (\zeta)  }{\zeta} =   \f{\vep}{2}  \li [ 1 -
\f{\nu \vep}{6 \si^4} + O (\vep^2) \ri ]  =   \f{\vep}{2}  \li [ 1 +
O (\vep) \ri ]
\]
as $\vep \to 0$. \eeT

See Appendix \ref{Aspstruecureapp} for a proof.

Applying Theorem \ref{simple} to independent bounded random
variables with known means, we have the following result.

\beT

\la{Hoeffdingmu}

Let $X_1, X_2,\cd$ be a sequence of independent random variables
such that $\Pr \{ 0 \leq X_i \leq 1 \} = 1$ and $\bb{E} [ X_i ] =
\mu_i$ for $i \in \bb{N}$.  Define $S_n = \sum_{i = 1}^n X_i$ and
$\ovl{\mu}_n = \f{ 1 }{n} \sum_{i = 1}^n \mu_i$ for $n \in \bb{N}$.
Define $\mcal{V}  (s, n) = n \ln ( \ovl{\mu}_n e^s + 1 - \ovl{\mu}_n
)$ for $s \in \bb{R}$ and $n \in \bb{N}$. Then, for all positive
integer $m$,
\[ \Pr \li \{  \sup_{n \in \bb{N} }
\li [ \zeta ( S_n - m \se ) - \mcal{V} ( \zeta, n) + \mcal{V} (
\zeta, m ) \ri ] \geq 0 \ri \} \leq \li [ \exp \li ( \se \ln
\f{\ovl{\mu}_m}{\se} + (1 - \se) \ln \f{ 1 - \ovl{\mu}_m }{1 - \se}
\ri ) \ri ]^m
\]
for all $\se \in (0, 1)$, where $\zeta =  \ln \f{\se (1 -
\ovl{\mu}_m)}{\ovl{\mu}_m (1 - \se)}$.

\eeT

See Appendix \ref{Hoeffdingmuapp} for a proof.

Applying Theorem \ref{simple} to independent bounded random
variables with known variances, we have the following result.

\beT

\la{Hoeffdingvar}

Let $X_1, X_2, \cd$ be independent random variables such that
$\bb{E} [ X_i ] = 0, \; \bb{E} [ | X_i |^2 ] \leq \si_i^2$,  and
$X_i \leq b$ for $i \in \bb{N}$.  Define $\nu_n = \f{1}{n} \sum_{i =
1}^n \si_i^2$ and $\mcal{V} (s, n) = n \ln \li [ \f{b^2}{b^2 +
\nu_n} \exp \li ( - \f{\nu_n}{b} s \ri ) + \f{\nu_n}{b^2 + \nu_n}
e^{b s} \ri ]$ for $n \in \bb{N}$. Then, \bee \Pr \li \{  \sup_{n
\in \bb{N} } \li [ \zeta ( S_n - m \vep ) - \mcal{V} ( \zeta, n) +
\mcal{V} ( \zeta, m ) \ri ] \geq 0 \ri \} \leq  \li [ \li ( 1 + \f{
b \vep }{ \nu_m } \ri )^{ - \f{ \nu_m + b \vep  }{ b^2 + \nu_m }   }
\; \li ( 1 - \f{\vep}{b} \ri )^{ - \f{ b^2 - b \vep }{ b^2 + \nu_m }
} \ri ]^m \eee for $ 0 < \vep < b$, where
\[
\zeta = \f{b}{b^2 + \nu_m} \ln \f{ 1 + \f{ \vep b}{\nu_m}  }{ 1 -
\f{\vep}{b} }.
\]
\eeT

See Appendix \ref{Hoeffdingvarapp} for a proof.

Applying Theorem \ref{simple} to independent random variables with
normal distributions, we have the following result.

\beT

\la{Umiformnormal}

 Let $X_1, X_2,\cd$ be a sequence of
independent random variables with normal distribution  such that
$\bb{E} [ X_i ] = \mu_i$ and $\bb{E} [ | X_i - \mu_i |^2 ] = \nu_i$
for $i \in \bb{N}$. Define $S_n = \sum_{i = 1}^n X_i$ and
$\ovl{\mu}_n = \f{ 1 }{n} \sum_{i = 1}^n \mu_i, \; \ovl{\nu}_n = \f{
1 }{n} \sum_{i = 1}^n \nu_i$ for $n \in \bb{N}$. Define $\mcal{V}
(s, n) = n \li ( \ovl{\mu}_n s + \f{ \ovl{\nu}_n  s^2}{2} \ri )$ for
$s \in \bb{R}$ and $n \in \bb{N}$. Then, for all positive integer
$m$,
\[ \Pr \li \{  \sup_{n \in \bb{N} }
\li [ \zeta ( S_n - m \se ) - \mcal{V} ( \zeta, n) + \mcal{V} (
\zeta, m ) \ri ] \geq 0 \ri \} \leq \li [ \exp \li ( - \f{  | \se -
\ovl{\mu}_m |^2 }{ 2 \ovl{\nu}_m} \ri )  \ri ]^m
\]
for all $\se \in \bb{R}$, where $\zeta =  \f{\se -
\ovl{\mu}_m}{\ovl{\nu}_m}$.

\eeT

See Appendix \ref{Umiformnormalapp} for a proof.

Applying Theorem \ref{simple} to independent Poisson random
variables, we have the following result.

\beT

\la{UniformPoisson}

Let $X_1, X_2,\cd$ be a sequence of independent Poisson random
variables such that  $\bb{E} [ X_i ] = \lm_i$ for $i \in \bb{N}$.
Define $S_n = \sum_{i = 1}^n X_i$ and $\ovl{\lm}_n = \f{ 1 }{n}
\sum_{i = 1}^n \lm_i$ for $n \in \bb{N}$. Define $\mcal{V}  (s, n) =
n \ovl{\lm}_n (e^s - 1)$ for $s \in \bb{R}$ and $n \in \bb{N}$.
Then, for all positive integer $m$,
\[ \Pr \li \{  \sup_{n \in \bb{N} }
\li [ \zeta ( S_n - m \se ) - \mcal{V} ( \zeta, n) + \mcal{V} (
\zeta, m ) \ri ] \geq 0 \ri \} \leq \li [ \exp \li ( \se -
\ovl{\lm}_m + \se \ln \f{\ovl{\lm}_m}{\se} \ri ) \ri ]^m
\]
for all $\se \in (0, \iy)$, where $\zeta =  \ln \f{\se}{\ovl{\lm}_m
}$.

\eeT

See Appendix \ref{UniformPoissonapp} for a proof.

\subsection{Using Moment Generating Functions}

Making use of moment generating functions pertained to vector
magnitude of random vectors, we have obtained the following results.

 \beT

Let $X, X_1, \cd, X_n$ be i.i.d. zero-mean random vectors. Let $Z$
be a zero-mean random variable in $\bb{R}$ such that $\bb{E} [ Z^k ]
\geq 1$ for $k \geq 2$. Assume that there exists a function
$\mscr{M} (s)$ such that $\bb{E} [ e^{ s Z || X ||} ] \leq \mscr{M}
(s)$ for all $s \in (- \tau, \tau)$, where $\tau
> 0$. Then, for any $\vep > 0$, \be \la{mainthm}
\Pr \li \{ \li | \li |  \f{1}{n} \sum_{i=1}^n X_i  \ri | \ri |  \geq
\vep \ri \} \leq \Pr \li \{ \max_{1 \leq \ell \leq n }  \li | \li |
\sum_{i = 1}^\ell X_i \ri | \ri | \geq n \vep \ri \} \leq \inf_{t
\in (0, \tau) } e^{ - n t \vep } \li \{ [ \mscr{M} (t) ]^n + [
\mscr{M} (- t) ]^n \ri \}.  \ee In particular, (\ref{mainthm}) holds
if the associated random variable $Z$ has a distribution such that
$\Pr \li \{ Z = \varphi \ri \} = \f{1}{\sq{5} \; \varphi}$ and $\Pr
\li \{ Z = - \f{1}{\varphi} \ri \} = \f{\varphi}{\sq{5}}$, where
$\varphi = \f{ 1 + \sq{5} }{2}$ is the golden ratio.

\eeT

In the case that the moment generating function of the magnitude of
a random vector exists, we have the following result.

\beT

Let $X, X_1, \cd, X_n$ be i.i.d. zero-mean random vectors such that
$\bb{E} [ e^{ s || X ||}  ] = g (s)$ for all $s \in (- \tau, \tau)$,
where $\tau > 0$. Let $\varphi = \f{ 1 + \sq{5} }{2}$ be the golden
ratio. Define \[ h (t, \vep, n) = e^{ - n t \vep } \li \{ \li [ \f{g
( \varphi t) }{ \varphi }  + \varphi g \li ( - \f{t}{ \varphi } \ri
) \ri ]^n + \li [ \f{g ( - \varphi t) }{ \varphi }  + \varphi g \li
( \f{t}{ \varphi } \ri ) \ri ]^n \ri \}
\]
for $\vep > 0$ and $t \in (0, \tau)$.  Then, for any $\vep > 0$,
\[
\Pr \li \{ \li | \li |  \f{1}{n} \sum_{i=1}^n X_i  \ri | \ri |  \geq
\vep \ri \} \leq  \Pr \li \{ \max_{1 \leq \ell \leq n }  \li | \li |
\sum_{i = 1}^\ell X_i \ri | \ri | \geq n \vep \ri \} \leq \f{1}{
\sq{5^n}} \times \inf_{t \in (0, \tau) }  h (t, \vep, n),
\]  where $h (t, \vep, n)$ is a convex function of $t \in (0, \tau)$ for fixed $\vep > 0$ and $n$.  \eeT

\subsection{Bounded Random Vectors}

Because of physical limitations, the magnitude of uncertainty
affecting systems are actually bounded.  Hence, it is of particular
importance to investigate the concentration phenomena of bounded
random vectors.

\subsubsection{Using Information of Support}

In the case that the bounds on the magnitude of random vectors are
available, we have the following result.

\beT

Let $X_1, \cd, X_n$ be independent zero-mean random vectors such
that $\Pr \{ || X_i || \leq r_i \} = 1$ for $i = 1, \cd, n$.
 Then, for all $\vep > 0$,
 \[
\Pr \li \{ \li | \li |  \f{1}{n} \sum_{i=1}^n X_i  \ri | \ri |  \geq
\vep \ri \} \leq \Pr \li \{ \max_{1 \leq \ell \leq n }  \li | \li |
\sum_{i = 1}^\ell  X_i \ri | \ri | \geq n \vep \ri \} \leq  2 \exp
\li ( - \f{ 2 n \vep^2 }{ 5 V} \ri ),
 \]
where $V = \f{1}{n} \sum_{i =1}^n r_i^2$.

 \eeT

If the diameters of the domain containing random vectors are known,
we have the following result.

 \beT
Let $X_1, \cd, X_n$ be independent zero-mean random vectors such
that $X_i$ has a support of diameter $D_i$ for $i = 1, \cd, n$.
Then, for all $\vep > 0$,
 \[
\Pr \li \{ \li | \li |  \f{1}{n} \sum_{i=1}^n X_i  \ri | \ri | \geq
\vep \ri \} \leq \Pr \li \{ \max_{1 \leq \ell \leq n }  \li | \li |
\sum_{i = 1}^\ell X_i \ri | \ri | \geq n \vep   \ri  \} \leq  2 \exp
\li ( - \f{ 2 n \vep^2 }{ 5 V} \ri ),
 \]
 where $V = \f{1}{n} \sum_{i=1}^n D_i^2$.

 \eeT

For vector-valued martingales of bounded increments, we have derived
maximal inequalities as follows.

\beT
 Suppose $\{ X_k: k = 0, 1, 2, 3, \cd \}$ is a vector-valued martingale and $\Pr \{ || X_k - X_{k-1} || \leq c_k \} = 1$
 for $k \in \bb{N}$.  Then,
 \[
\Pr \li \{  || X_n - X_0 || \geq \vep  \ri \} \leq 2 \exp \li ( -
\f{ 2 \vep^2}{ 5 \sum_{k = 1}^n c_k^2 } \ri )
 \]
 for all positive integers $n$ and all positive reals $\vep$.

 \eeT

 \subsubsection{Using Information of Support and Variance}

To make use of the information of each component of random vectors,
we have the following results.

\beT
 Let $X = [\bs{x}_1, \cd, \bs{x}_d ]$ be a zero-mean random vector such that $\bb{E} [ || X  ||^2 ] \leq \si^2$,
 the components $\bs{x}_1, \cd, \bs{x}_d$ are mutually independent,   and
 $\Pr \{ | \bs{x}_i | \leq r_i \} = 1$ for $i = 1, \cd, d$.  Then,
 \[
\Pr \{  || X || \geq \vep \} \leq  \exp \li (   - \f{ 2 | \vep^2 -
\si^2 |^2  }{  \sum_{i=1}^d r_i^4 }   \ri )
 \]
for $\vep > \si$.

\eeT

If we know the range of each component of random vectors, we have
the following result.

\beT
 Let $X = [\bs{x}_1, \cd, \bs{x}_d ]$ be a zero-mean random vector such that the components $\bs{x}_1, \cd, \bs{x}_d$ are mutually
 independent and that $\Pr \{ a_i \leq \bs{x}_i \leq b_i \} = 1$ for $i = 1, \cd,
 d$.  Define $\si^2 = \sum_{i = 1}^d | a_i b_i | $.
Then,
 \[
\Pr \{  || X - \mu || \geq \vep \} \leq  \exp \li (   - \f{ 2 |
\vep^2 - \si^2 |^2  }{  \sum_{i=1}^d | b_i - a_i |^4 }   \ri )
 \]
for $\vep > \si$.

\eeT

Making use of the variance information of random vectors, we have
derived simple exponential inequalities as follows.

\beT

\la{expbound3388} Let $X_1, X_2, \cd$ be independent zero-mean
random vectors  such that for $n \in \bb{N}$, \bee  \sum_{i = 1}^n
\bb{E} [ || X_i ||^2 ] \leq s_n^2, \qqu  \Pr \{ || X_i || \leq c_n
s_n \; \tx{for} \; i = 1, \cd, n \} = 1,  \eee where $c_n > 0$ and
$s_n
> 0$.  Let $\varphi = \f{ 1 + \sq{5} }{2}$ be the golden ratio.
Then,
\[
\Pr \li \{ \max_{1 \leq \ell \leq n }  \li | \li | \sum_{i = 1}^\ell
X_i \ri | \ri | \geq x s_n \ri \} \leq 2 \exp \li ( - \f{x^2}{2} \li
( 1 - \f{x \varphi c_n}{2} \ri ) \ri )
\]
for $0 < x < \f{1}{\varphi c_n}$.

 \eeT

Making use of the variance and range information of random vectors,
we have derived tight inequalities as follows.

\beT

\la{gen8899}

 Let $X_1, \cd, X_n$ be independent zero-mean random vectors such
that $\sum_{i =1}^n \bb{E} [  || X_i ||^2 ] \leq n \si^2$ and $\Pr
\{ || X_i || \leq r \} = 1$ for $i = 1, \cd, n$, where $\si \geq 0$
and $r > 0$. Let $\varphi = \f{1 + \sq{5}}{2}$ be the golden ratio.
Then, $\Pr \{ \li | \li | \f{1}{n} \sum_{i=1}^n X_i  \ri | \ri |  >
r \} = 0$ and

{\small $\Pr \li \{ \li | \li |  \f{1}{n} \sum_{i=1}^n X_i \ri | \ri
| \geq \vep \ri \}  \leq  \Pr \li \{ \max_{1 \leq \ell \leq n } \li
| \li | \sum_{i = 1}^\ell  X_i \ri | \ri | \geq n \vep   \ri \}$}

{\small $= \inf_{t > 0} \; e^{ - n t \vep } \li \{ \li [ \f{
(\varphi r)^2 }{ \si^2 + (\varphi r)^2} \exp \li ( - \f{ t \si^2
}{\varphi r} \ri )
 + \f{ \si^2 }{\si^2 + (\varphi r)^2} \exp ( t \varphi r )
\ri ]^n  + \li [ \f{ r^2 }{r^2 + (\varphi \si)^2 } \exp \li ( - \f{
t \varphi \si^2 }{r} \ri ) + \f{ (\varphi \si)^2 }{r^2 + (\varphi
\si)^2} \exp \li ( \f{ t r}{\varphi} \ri ) \ri ]^n \ri \}$}

{\small $\leq   2 \li [ \li (  \f{  \si^2 }{  \si^2 + \varphi r \vep
} \ri )^{ \si^2 + \varphi r \vep } \li ( 1 - \f{\vep}{\varphi r} \ri
)^{ \varphi r \vep - (\varphi r)^2} \ri ]^{\f{n}{\si^2 + (\varphi
r)^2 }} \leq  2 \li [ \li (  \f{ \si^2 }{ \si^2 + \varphi r \vep }
\ri )^{ \si^2 + \varphi r \vep  } \exp \li ( \varphi r \vep \ri )
\ri ]^{ \f{n}{(\varphi r)^2} } \leq  2 \exp \li (  - \f{  n \vep^2
}{ 2 ( \si^2 + \f{\varphi r \vep}{3} ) } \ri )$ } for $0 < \vep \leq
r$.

 \eeT

To apply Theorem \ref{gen8899}, we need to bound $|| X - \mu ||$ and
$\bb{E} [ || X - \mu ||^2 ]$.  For this purpose, we have the
following result.

\beT

Let $X$ be a random vector with mean $\mu = \bb{E} [ X ]$ and a
support of diameter $D$. Then, $|| X - \mu || \leq D$ and $\bb{E} [
|| X - \mu ||^2 ] \leq \f{D^2}{2}$. \eeT

If random vector $X$ is bounded within an ellipse, we have the
following result.

 \beT

 \la{extenBD}

Let $X$ be a random vector such that $|| A X + b || \leq c$, where
$A$ is an invertible matrix. Then,
 \[
|| X - \mu || \leq  || A^{-1} || \times [ c + || A  \mu + b || ],
\qqu  \bb{E} [  || X - \mu ||^2  ] \leq || A^{-1} || \times [ c^2 -
|| A \mu  + b ||^2],
 \]
where $\mu = \bb{E} [ X ]$.

 \eeT

It should be noted that Theorem \ref{extenBD} is an extension of
Bhatia-Davis inequality  \cite{Bhatia}.

\section{Stability of Uncertain Dynamic Systems}

\begin{figure}[htbp]
\centerline{\psfig{figure=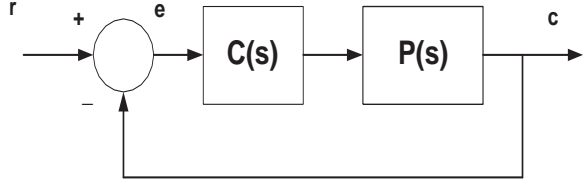, height=1.1in, width=3.5in }}
\caption{Uncertain System } \label{fig_08}
\end{figure}

In this section, we shall apply the proposed theory of inferencing
function of uncertainties to study the stability of uncertain
systems. Consider a system which has been studied in \cite{GS} by a
deterministic approach. The system is shown in Figure~\ref{fig_08}.

The compensator is $C(s)=\frac{s+2}{s+10}$ and the plant is {\small
$P(s)=\frac{800(1+0.1\eta_1)}{s(s+4+0.2\eta_2)(s+6+0.3\eta_3)}$}
with parametric uncertainty $| \eta_i | \leq 0.16, \qu | \bb{E} [
\eta_i ] | < 0.05$ for $i = 1, \; 2, \; 3$.  The transfer function
of the system is {\small $T(s) = \f{  C (s) P (s) }{ 1 + C (s) P (s)
}$}. The characteristic polynomial of the system is \bee  s (s + 10)
(s+4+0.2 \eta_2)(s+6+0.3\eta_3) + 800(1+0.1\eta_1) (s + 2) = s^4 +
a_1 s^3 + a_2 s^2 + a_3 s + a_4,  \eee where \bee  &  &  a_1 = 20 +
0.2 \eta_2 +  0.3\eta_3,   \qqu \qqu \qqu  a_2 = (4+0.2 \eta_2)
( 6+0.3\eta_3  ) +  10 (10 + 0.2 \eta_2 +  0.3\eta_3) , \\
&  & a_3 = 10 (4+0.2 \eta_2) ( 6+0.3\eta_3  )  + 800(1+0.1\eta_1),
\qqu \qqu  a_4 = 1600(1+0.1\eta_1). \eee By the Routh stability
criterion, the system is stable if and only if
\[
a_1 > 0, \qqu a_1 a_2 - a_3 > 0, \qqu (a_1 a_2 - a_3 ) a_3 - a_1^2
a_4 > 0, \qqu a_4 > 0,
\]
that is, $h(\eta_1, \eta_2, \eta_3) > 0$, where $h(\eta_1, \eta_2,
\eta_3) = \min \{ a_1, \; a_4, \;  a_1 a_2 - a_3, \; (a_1 a_2 - a_3
) a_3 - a_1^2 a_4 \}$.  Hence, if we define
\[
X = [ \eta_1, \; \eta_2, \; \eta_3], \qqu \bs{f} (X) = X, \qqu g(X)
= \bb{I}_{ \{h(\eta_1, \eta_2, \eta_3) > 0 \} },
\]
\[
\mscr{A} = \{ (x_1, x_2, x_3): |x_i| < 0.16, \; i = 1, \; 2, \; 3
\}, \qqu \mscr{B} = \{ (x_1, x_2, x_3): |x_i| < 0.05, \; i = 1, \;
2, \; 3 \},
\]
\[
\mscr{C} = \{ (x_1, x_2, x_3): \; h (x_1, \; x_2, \; x_3 ) \leq 0
\},
\]
then
\[
\Pr \{ \tx{The system is unstable} \} = \Pr \{ X \in \mscr{C} \},
\]
subject to
\[
\Pr \{ X \in \mscr{A} \} = 1, \qqu \bb{E} [ \bs{f} ( X) ] \in
\mscr{B}.
\]
Therefore, we can apply Theorem \ref{VIPProb} to compute a
deterministic bound for $\Pr \{ \tx{The system is unstable} \}$.
With less than $0.05$ second, we obtained such upper bound as
$0.00031$ by a computer program which implements linear programming
embedded with the gradient search and the branch and bound
algorithms.

\section{Conclusion}

In this paper, we have developed a general theory for inferring
uncertainty.  We have applied the general theory to investigate
concentration phenomena of random vectors. Uniform exponential
inequalities and multidimensional probabilistic inequalities have
been developed which can be useful for the analysis of control and
decision affected by uncertainty. We have derived computable tight
bounds for the  expected values of functions of uncertainty which
represent performance of systems. The applications of such results
are illustrated by an investigation of the stability of an uncertain
system.

\appendix

\section{Proof of Theorem \ref{theorem88388chen} } \la{theorem88388chenapp}

Note that since all elements in $\mscr{Y}$ are discrete random
vectors, the associated expectation $\bb{E} [ g (Y)]$ of any $Y \in
\mscr{Y}$  must exist. Hence, $\sup_{ Y \in \mscr{Y} } \bb{E} [ g
(Y) ]$ is well-defined provided that $\mscr{Y}$ has at least one
element. Therefore, it suffices to show that the family $\mscr{Y}$
contains at least one element $Y$ with $\bb{E} [ g (Y) ] \geq \bb{E}
[ g (X) ]$.  Define $S = \{ (u, v):  u = \bs{f} (x), \; v = g (x),
\; x \in \mscr{A} \}$. Then, $\Pr \{ (\bs{f} (X), \; g (X) ) \in S
\} = 1$. Note that the convex hull of $S$, denoted by $\mrm{conv}
(S)$, is convex. By assumption, both $\bb{E} [ \bs{f} (X) ]$ and
$\bb{E} [ g (X) ]$ exist. Hence, by Theorem \ref{fundamental},
\[
( \bb{E} [ \bs{f} (X) ], \; \bb{E} [ g (X) ] ) \in \mrm{conv} (S).
\]
Note that $S$ is a subset of $(k+1)$-dimensional vector space.
According to Carathéodory's theorem, there exists $m \leq k + 2$
points, $x_1, \cd, x_m$ in $S$ such that $( \bb{E} [ \bs{f} (X) ],
\; \bb{E} [ g (X) ] )$ is a convex combination of $x_1, \cd, x_m$.
The points $x_1, \cd, x_m$ are vertexes of the simplex which
consists of all convex combinations of $x_1, \cd, x_m$.  Consider
half-line $\{ (u, v): u = \bb{E} [ \bs{f} (X) ], \; v \geq \bb{E} [
g (X) ] \}$.  There must exist $w \geq \bb{E} [ g (X) ]$ such that
$(\bb{E} [ \bs{f} (X) ], w)$ lie in a proper face of the simplex.

Without loss of generality, let $x_1, \cd, x_{m-1}$ be the vertex of
such proper face. Then, there exist nonnegative numbers $p_1, \cd,
p_{m-1}$ such that $\sum_{i = 1}^{m -1} p_i = 1$ and that
\[
\bb{E} [ \bs{f} (X) ] = \sum_{i =1}^{m -1} p_i \bs{f} (x_i), \qqu w
= \sum_{i =1}^{m -1} p_i g (x_i).
\]
Hence, we can define a discrete random vector $Y$ of $(m-1) \leq k +
1$ possible values such that $\Pr \{ Y = x_i \} = p_i$ for $i = 1,
\cd, m -1$.  Clearly,
\[
\Pr \{ Y \in \mscr{A} \} = 1, \qqu \bb{E} [ \bs{f} (Y) ]  = \sum_{i
=1}^{m -1} p_i \bs{f} (x_i) = \bb{E} [ \bs{f} (X) ] \in \mscr{B},
\qqu \bb{E} [ g (Y) ] = \sum_{i =1}^{m -1} p_i g (x_i) = w \geq
\bb{E} [ g (X) ].
\]
This shows that the family $\mscr{Y}$ contains at least one element
$Y$ with $\bb{E} [ g (Y) ] \geq \bb{E} [ g (X) ]$.  The proof of the
theorem is thus complete.

\section{Proof of Theorem \ref{VIPProb}} \la{VIPProbapp}

For $y \in \mscr{A}$, define $g(y)$ such that $g(y) = 1$ if $y \in
\mscr{C}$ and that $g(y) = 0$ if $y \in \mscr{A} \setminus
\mscr{C}$. According to Theorem \ref{theorem883}, we have {\small
\bee \sup_{ X \in \mscr{X} }  \Pr \{ X \in \mscr{C}  \} & = & \sup
\li \{ \sum_{\ell = 1}^{k + 1} \se_\ell g (y_\ell) : \; \se_\ell
\geq 0 \; \tx{and} \; y_\ell \in \mscr{A} \; \tx{for} \; \ell = 1,
\cd, k + 1, \qu \sum_{\ell = 1 }^{k + 1} \se_\ell = 1, \qu
\sum_{\ell = 1}^{k + 1} \se_\ell \bs{f} ( y_\ell )
\in \mscr{B} \ri \}\\
& = &  \max \{ Q_i: 1 \leq i \leq k + 1  \}, \eee} where {\small
\[
Q_i =  \sup \li \{ \sum_{\ell = 1}^{k + 1} \se_\ell g (y_\ell) : \;
\se_\ell \geq 0 \; \tx{and} \; y_\ell \in \mscr{A} \; \tx{for} \;
\ell = 1, \cd, k + 1, \qu \sum_{\ell = 1}^{k + 1} \se_\ell = 1, \qu
\sum_{\ell = 1}^{k + 1} g (y_\ell) = i, \qu \sum_{\ell = 1}^{k + 1}
\se_\ell \bs{f} ( y_\ell ) \in \mscr{B} \ri \} \] } for $i=1, \cd,
k+1$. Define $E_i = \{ (b_1, \cd, b_{k+1}): \sum_{\ell = 1}^{k+1}
b_\ell = i, \; \tx{where} \; b_\ell \in \{ 0, 1 \} \; \tx{for} \;
\ell = 1, \cd, k +1 \}$ for $i = 1, \cd, k + 1$. Then,  $Q_i  = \max
\{ h (b_1, \cd, b_{k+1}):  (b_1, \cd, b_{k+1}) \in E_i  \}$, where
$h (b_1, \cd, b_{k+1})$ is defined as {\small $\sup \{ \sum_{\ell =
1}^{k + 1} \se_\ell b_\ell : \; \se_\ell \geq 0 \; \tx{and} \;
y_\ell \in \mscr{A},  \; g(y_\ell) = b_\ell \; \tx{for} \; \ell = 1,
\cd, k + 1, \qu \sum_{\ell = 1}^{k + 1} \se_\ell = 1, \qu \sum_{\ell
= 1}^{k + 1} \se_\ell \bs{f} ( y_\ell ) \in \mscr{B} \}$,} for $i=1,
\cd, k+1$. Consider $(b_1, \cd, b_{k+1}) \in E_i$ such that
$b_{\ell_t} = 1$ for $1 \leq t \leq i$ and $b_{\ell_t} = 0$ for $i <
t \leq k + 1$. Define $x_t = y_{\ell_t}$ and $\vse_t = \se_{\ell_t}$
for $t = 1, \cd, k + 1$. Then, $h (b_1, \cd, b_{k+1}) $ is equal to
{\small \bee &  & \sup \li \{ \sum_{\ell = 1}^i \vse_\ell: \;
\vse_\ell \geq 0  \; \tx{for} \; 1 \leq \ell \leq k + 1, \qu  x_\ell
\in \mscr{C} \; \tx{for} \; 1 \leq \ell \leq i,
 \qu  x_\ell \in \mscr{A} \setminus \mscr{C} \; \tx{for} \; i < \ell \leq k + 1, \ri. \\
&  &  \qu \qu \qu \qu \qu \qu \li.  \sum_{\ell = 1 }^{k + 1}
\vse_\ell = 1, \qu \sum_{\ell = 1}^{k + 1} \vse_\ell \bs{f} ( x_\ell
) \in \mscr{B} \ri \}, \eee} which is the same as $P_i$.   Hence, we
have established that $h (b_1, \cd, b_{k+1}) = P_i$ holds for all
$(b_1, \cd, b_{k+1}) \in E_i$ for $1 \leq i \leq k + 1$.  It follows
that $Q_i = P_i$ for $1 \leq i \leq k + 1$. Therefore,  $\sup_{ X
\in \mscr{X} } \Pr \{ X \in \mscr{C}  \} = \max \{ Q_i: 1 \leq i
\leq k + 1  \} = \max \{ P_i: 1 \leq i \leq k + 1  \}$.  This
completes the proof of the theorem.

\section{Proof of Theorem \ref{Them9}} \la{Them9_app}

Define $W_t = \exp ( s (X_t - X_0) - \varphi (s) \mcal{V}_t )$ for $
t \geq 0$ and $s \in (0, b)$. Then, for all $s \in (0, b)$ and
arbitrary $t^\prime \geq t \geq 0$, we have {\small \bee &  & \bb{E}
[ W_{t^\prime} \mid \mscr{F}_t] = \bb{E} \li [  \exp ( s (
X_{t^\prime} - X_0 ) - \varphi (s) \mcal{V}_{t^\prime} ) \mid
\mscr{F}_t \ri ]
 = \bb{E} \li [  \exp ( s (X_{t^\prime} - X_t) - \varphi (s) (\mcal{V}_{t^\prime} - \mcal{V}_t)  ) \; W_t \mid \mscr{F}_t \ri ]\\
&  & =  W_t  \exp ( - \varphi (s) (\mcal{V}_{t^\prime} - \mcal{V}_t)
) \; \bb{E} \li [  \exp ( s (X_{t^\prime} - X_t)  ) \mid \mscr{F}_t
\ri ] \leq W_t. \eee} Hence, for any $s \in (0, b)$, $( W_t,
\mscr{F}_t )_{t \in \bb{R}^+ }$ is a super-martingale with $\bb{E}
[W_0] = \bb{E} [ \exp ( - \varphi (s) \mcal{V}_0 ) ] \leq 1$.  By
the assumption on the continuity of the sample paths of $\{ s (X_t -
X_0) - \varphi (s) \mcal{V}_t \}_{t > 0}$, we have that almost all
sample paths of $( W_t )_{t \in \bb{R}^+}$ is right-continuous.

To prove (\ref{genb12}), note that for any $s \in (0, b)$ and real
number $\ga > 0$, {\small \bel  &  & \Pr \li \{ \sup_{t > 0} \li [
X_t - X_0 - \ga \mcal{V}_\tau - \f{ \varphi(s) }{s} (\mcal{V}_t -
\mcal{V}_\tau) \ri ] \geq 0 \ri \}   =
\Pr \li \{ \sup_{t > 0} \li [ X_t - X_0 - \ga \mcal{V}_\tau - \f{ \varphi(s) }{s} (\mcal{V}_t - \mcal{V}_\tau) \ri ] s \geq 0 \ri \} \nonumber\\
&  & = \Pr \li \{ \sup_{t > 0} \li [ s ( X_t - X_0 ) - \varphi(s) \mcal{V}_t -
\ga s \mcal{V}_\tau  + \varphi(s) \mcal{V}_\tau \ri ] \geq 0 \ri \}  =
\Pr \li \{ \sup_{t > 0} \li [ s ( X_t - X_0 ) - \varphi(s) \mcal{V}_t \ri ] \geq \ga s \mcal{V}_\tau  - \varphi(s) \mcal{V}_\tau  \ri \}\nonumber\\
&  & = \Pr \li \{ \sup_{t > 0} W_t  \geq \exp \li ( \ga s \mcal{V}_\tau  - \varphi(s) \mcal{V}_\tau \ri ) \ri \} \la{defm336}\\
&  & \leq \exp \li ( \varphi(s) \mcal{V}_\tau - \ga s \mcal{V}_\tau  \ri ) \la{marineq339}\\
&  &  = \li [ \exp \li ( \varphi(s) - \ga s  \ri ) \ri
]^{\mcal{V}_\tau}. \nonumber \eel}  Here, we have used the
definition of $W_t$ in (\ref{defm336}). The inequality
(\ref{marineq339}) follows from the super-martingale inequality.
This proves (\ref{genb12}) and thus (\ref{discoverb12}) immediately
follows.  This concludes the proof of the theorem.

\section{Proof of Theorem \ref{Them98} } \la{Them98app}

To prove Theorem \ref{Them98}, note that for all $s \in \mscr{S}$,
{\small \bee \Pr \li \{ \sup_{t > 0}  \li [ s ( X_t - \ga ) -
\mcal{V} (s, t) + g (s) \ri ] \geq 0  \ri \} & = & \Pr  \li \{
\sup_{t > 0} \li [ s X_t - \mcal{V} (s, t) \ri ] \geq \ga s - g (s) \ri \}\\
& = &  \Pr \li \{ \sup_{t > 0} \exp \li ( s X_t - \mcal{V} (s, t)
\ri ) \geq \exp (\ga s - g (s) ) \ri \}\\
&  \leq & \Pr \li \{ \sup_{t
> 0} \mcal{Z} (s, t)  \geq \exp ( \ga s - g (s) ) \ri \}. \eee} By the
supermartingale inequality, we have {\small \[ \Pr \li \{ \sup_{t >
0}  \li [ s ( X_t - \ga ) - \mcal{V} (s, t) + g (s) \ri ] \geq 0 \ri
\}  \leq \f{ \bb{E} [ \mcal{Z} (s, 0) ] }{ \exp ( \ga s - g (s) ) }
\leq \f{ 1 }{ \exp ( \ga s - g (s) ) } =  \exp (g (s) - \ga s )
\] }
for all $s \in \mscr{S}$. This proves (\ref{gen8833}), from which
the particular assertions  immediately follow.

\section{Proof of Theorem \ref{simple} } \la{simpleapp}

Define \[ \phi (s, n) = \ln \bb{E} [ e^{s Y_n} ], \qqu  \varphi (s,
n) = \ln \bb{E} [ e^{s X_n} ], \qqu \mcal{Z} (s, n) = \exp \li (  s
X_n -  \varphi (s, n) \ri )
\]
for $s \in \mscr{S}$ and $n \in \bb{N}$.  Clearly, $\mcal{Z} (s, 0)
= 1$ for all $s \in \mscr{S}$.  For $n \in \bb{N}$, let $\mscr{F}_n$
denote the $\si$-algrbra generated by $Y_1, \cd, Y_n$. Note that
\bee \bb{E} [ \mcal{Z} (s, n + 1) \mid \mscr{F}_n ] & = & \bb{E} [
\exp \li ( s X_n + s Y_{n+1} - \varphi (s, n) - \phi (s, n + 1) \ri
)
\mid \mscr{F}_n ] \\
& = & \bb{E} [ \exp \li (  s Y_{n+1} -  \phi (s, n + 1) \ri ) ]
\times \bb{E} [  \exp \li (  s X_n  -  \varphi (s, n) \ri ) \mid \mscr{F}_n ] \\
& = &  \mcal{Z} (s, n) \eee almost surely.  Hence,  for each $s \in
\mscr{S}$,  $\{ \mcal{Z} (s, n), \; n \in \bb{N} \}$ is a
martingale.   By assumption,
\[
\exp \li (  s X_n -  \mcal{V} (s, n) \ri ) \leq \exp \li (  s X_n -
\varphi (s, n) \ri ) = \mcal{Z} (s, n)
\]
for $s \in \mscr{S}$ and $n \in \bb{N}$.  Invoking assertion (II) of
Theorem \ref{Them98} yields the conclusion of the theorem.

\section{Proof of Theorem  \ref{Aspstruecure88689}} \la{Aspstruecure88689app}

We need some preliminary results.

\beL \la{Order3388689}

Define $h (s, \vep) =  \varphi (s) - s \vep$ for $s \in (a, b)$ and
$\lim_{ s \downarrow a } \f{ \varphi
 (s)} {s} < \vep < \lim_{ s \uparrow b }  \f{ \varphi
 (s)} {s}$.  Then, there exists a number $c$ such that $0 < c < \min \{
 |a|, b \}$ and that for $\vep$ satisfying $- \f{\varphi (-c)}{c} <
 \vep < \f{\varphi (c)}{c}$,   the infimum of $h (s, \vep)$ with
 respect to $s \in (a, b)$ is attained at the unique root,  $\zeta = \zeta (\vep)$, of the
 equation $\varphi^\prime (s) = \vep$ with respect to $s \in (-c, c)$.
 Moreover, $\zeta = \zeta (\vep) = \f{\vep}{ \al } + O(\vep^2)$ for $- \f{\varphi (-c)}{c} <
 \vep < \f{\varphi (c)}{c}$.

\eeL

 \bpf  Let
$\varphi^\prime(s), \; \varphi^{\prime \prime} (s)$, and
$\varphi^{(3)} (s)$ denote the first, second, and third derivatives
of $\varphi(s)$, respectively.  By assumption,
\[
\varphi (0) = 0, \qqu \varphi^\prime (0) = 0, \qqu \varphi^{\prime
\prime} (0) = \al > 0.
\]
Clearly, $h(s, \vep)$ is convex with respect to $s \in (a, b)$. By
the assumption on $\vep$ and the convexity of $\varphi (s)$, we have
 \[
\lim_{ s \downarrow a } \varphi^\prime (s) \leq \lim_{ s \downarrow
a } \f{ \varphi
 (s)} {s} < \vep < \lim_{ s \uparrow b }  \f{ \varphi
 (s)} {s} \leq \lim_{ s \uparrow b } \varphi^\prime (s).
 \]
Since $\varphi (s)$ is convex for $s \in (a, b)$, it follows that
$\varphi^\prime (s)$ is increasing for $s \in (a, b)$.  Hence, the
equation $\varphi^\prime (s) = \vep$ with respect to $s \in (a, b)$
has a unique root,  $\zeta = \zeta (\vep)$, for $\vep$ such that
$\lim_{ s \downarrow a } \f{ \varphi
 (s)} {s} < \vep < \lim_{ s \uparrow b }  \f{ \varphi
 (s)} {s}$.  Moreover,
\[
\zeta =  \bec \zeta (\vep)   < 0 & \tx{for} \; \lim_{ s \downarrow a
} \f{ \varphi (s)} {s} < \vep < 0,\\
 \zeta (\vep)  = 0 & \tx{for} \; \vep = 0,\\
 \zeta (\vep)  > 0 & \tx{for} \; 0 < \vep < \lim_{ s \uparrow b } \f{ \varphi (s)} {s}
 \eec
\]
Since $h(s, \vep)$ is convex with respect to $s \in (a, b)$, it
follows that
 \[
\inf_{s \in (a, b)} h (s, \vep) = h (\zeta, \vep).
 \]
By assumption,  $\varphi^{\prime \prime} (0) = \al > 0$.  Since
$\varphi^{\prime \prime} (s)$ is a continuous function of $s \in (a,
b)$, it follows that there exists a positive number $c$ such that $0
< c < \min \{ |a|, \; b \}$, the third derivative $\varphi^{(3)}
(s)$ is continuous, and that $\varphi^{\prime \prime} (s) \geq
\f{\al}{2}$ for all $s \in [ - c, c ]$.  By Taylor series expansion
formula,
\[
\varphi^\prime (s)  = \varphi^\prime (0) + \varphi^{\prime \prime}
(\eta s) s = \varphi^{\prime \prime} (\eta s) s \geq \f{\al}{2} s
\qu \tx{for $s \in (0, c]$,}
\]
where $\eta \in (0, 1)$ is a number dependent on $s$.
 Similarly, $\varphi^\prime (s) \leq \f{\al}{2} s$
for $s \in [-c, 0)$.   As a consequence of the convexity of $\varphi
(s)$, it must be true that $\varphi^\prime (- c) < \vep <
\varphi^\prime (c)$ for $\vep$ satisfying $- \f{\varphi (-c)}{c} <
 \vep < \f{\varphi (c)}{c}$.   Since $\varphi^\prime (s)$ is monotonically increasing with respect to
$s \in [- c, c ]$, it follows that $ - c < \zeta = \zeta (\vep) < c$
for $- \f{\varphi (-c)}{c} < \vep < \f{\varphi (c)}{c}$.   From now
on, we restrict $\vep$ to satisfy the
 constraint $- \f{\varphi (-c)}{c} <
 \vep < \f{\varphi (c)}{c}$.  Therefore,
\[
\vep = \varphi^\prime (\zeta) \geq \f{\al}{2} \zeta \qu \tx{for} \qu
0 < \vep < \f{\varphi (c)}{c}
\]
and \[ \vep = \varphi^\prime (\zeta) \leq \f{\al}{2} \zeta \qu
\tx{for} \qu - \f{\varphi (-c)}{c} < \vep < 0.
\]
Hence, $\li | \f{ \zeta }{\vep} \ri | = \f{ \zeta }{\vep} \leq \f{ 2
}{ \al }$ for $0 < \vep < \f{\varphi (c)}{c}$ and $- \f{\varphi
(-c)}{c} < \vep < 0$. This shows that \be \la{result88a689}
 \zeta = \zeta (\vep) =
O (\vep) \ee for $- \f{\varphi (-c)}{c} < \vep < \f{\varphi
(c)}{c}$. By Taylor series expansion formula,
\[
\varphi^\prime (s) =  \al s  + \f{1}{2} \varphi^{(3)} (\eta s) s^2
\qu \tx{for $s \in [- c, c]$,}
\]
where $\eta \in (0, 1)$.  Hence,
\[
\varphi^\prime (\zeta) =  \al \zeta + \f{1}{2} \varphi^{ (3)} (\eta
\zeta) \zeta^2
\]
for $- \f{\varphi (-c)}{c} < \vep < \f{\varphi (c)}{c}$,  where
$\eta \in (0, 1)$ is dependent on $\vep$.   Since $\varphi^{ (3) }
(s)$ is continuous with respect to $s \in [- c, c]$, there exists $K
> 0$ such that $| \varphi^{ (3) } (s) | \leq K$ for all $s \in [ - c, c
]$. Recall that  $h(s, \vep)$ is minimized at $s = \zeta = \zeta
(\vep)$ such that $\varphi^\prime (\zeta) = \vep$.   Hence, $\zeta$
satisfies the equation
\[
\al \zeta + \f{1}{2} \varphi^{(3)} (\eta \zeta) \zeta^2 = \vep
\]
and thus \be \la{result88b689}
 \zeta = \f{\vep}{ \al}  - \f{1}{2 \al} \varphi^{(3)} (\eta \zeta) \zeta^2 \ee
for $- \f{\varphi (-c)}{c} < \vep < \f{\varphi (c)}{c}$. Since $|
\varphi^{(3)} (\eta \zeta) | \leq K$, it follows from
(\ref{result88a689}) and (\ref{result88b689}) that \[ \zeta  =
\f{\vep}{ \al } - \f{1}{2 \al} \varphi^{(3)} (\eta \zeta) [ O (\vep)
]^2  =  \f{\vep}{ \al } + O (\vep^2)  \] for $\vep$ satisfying $-
\f{\varphi (-c)}{c} < \vep < \f{\varphi (c)}{c}$. This completes the
proof of the lemma.

\epf

\bsk

We are now in a position to prove the theorem.  The inequality
(\ref{firstinequality689}) immediately follows from assertion (II)
of Theorem \ref{Them98}. The inequality (\ref{secondinequality689})
follows from (\ref{firstinequality689}) and the convexity of
$\varphi(s)$. It remains to investigate the asymptotic expression of
the probability bound $\li [ \exp  ( \varphi (\zeta) - \vep \zeta )
\ri ]^{\mcal{V}_\tau}$  and the ratio $\f{\varphi(\zeta)}{\zeta}$.

Let $c$ be the number in the context of Lemma \ref{Order3388689}.
Using Taylor series expansion formula, we have
\[ \varphi (s) = \varphi (0) + \varphi^\prime (0) s + \f{1}{2} \varphi^{\prime
\prime} (0) s^2 +  O (s^3) = \f{\al}{2} s^2 +  O (s^3) \qu \tx{for
$s \in [- c, c]$. }
\]
Making use of this expression of $\varphi (s)$ and the expression of
$\zeta$ in Lemma \ref{Order3388689}, we have \[  \varphi (\zeta)  =
 \f{\al}{2} \zeta^2 + O (\zeta^3)\\
  =  \f{\al}{2} \li [ \f{\vep}{ \al }
+ O(\vep^2)  \ri ]^2  +  O( \vep^3 ) \\
 =   \f{\vep^2}{2 \al}  + O (\vep^3), \] \[ h (\zeta, \vep)
=  - \vep \zeta + \varphi (\zeta)  =   - \f{\vep^2}{ \al } +
O(\vep^3) + \f{\vep^2}{2 \al} + O (\vep^3)  =   - \f{\vep^2}{2 \al}
+  O (\vep^3), \] \bee \li [   \exp  ( \varphi (\zeta) - \vep \zeta
) \ri ]^{\mcal{V}_\tau} & = & \li [ \exp
( h (\zeta, \vep)  ) \ri ]^{\mcal{V}_\tau}\\
& = & \li [  \exp \li (  - \f{\vep^2}{2 \al} +  O (\vep^3) \ri ) \ri ]^{\mcal{V}_\tau}\\
& = & \li [  \exp \li (  - \f{\vep^2}{2 \al}  \ri )  \exp \li (  O
(\vep^3) \ri )\ri ]^{\mcal{V}_\tau}\\
& = & \li [  \exp \li (  - \f{\vep^2}{2 \al} \ri )  \ri
]^{\mcal{V}_\tau} \li [ \exp \li (  O
(\vep^3) \ri )\ri ]^{\mcal{V}_\tau}\\
& = & \li [  \exp \li (  - \f{\vep^2}{2 \al}  \ri )  \ri ]^{\mcal{V}_\tau} \li [  1 +  O (\vep^3)  \ri ]^{\mcal{V}_\tau}\\
& = & \li [  1 +  O (\vep^3)  \ri ] \li [  \exp \li (  -
\f{\vep^2}{2 \al} \ri )  \ri ]^{\mcal{V}_\tau},
 \eee
  and  \bee
\f{\varphi (\zeta)}{\zeta}  =  \f{ \f{\vep^2}{2 \al}  + O (\vep^3)
}{ \f{\vep}{ \al } + O(\vep^2) }  =   \f{\vep}{2}  \li [ 1 +  O
(\vep) \ri ] \eee for $\vep$ satisfying $- \f{\varphi (-c)}{c} <
\vep < \f{\varphi (c)}{c}$. This completes the proof of the theorem.

\section{Proof of Theorem  \ref{Aspstruecure}} \la{Aspstruecureapp}

We need some preliminary results.

\beL \la{Order3388}

Define $h (s, \vep) =  \varphi (s) - s \vep$ for $s \in (a, b)$ and
$\lim_{ s \downarrow a } \f{ \varphi
 (s)} {s} < \vep < \lim_{ s \uparrow b }  \f{ \varphi
 (s)} {s}$.  Then, there exists a number $c$ such that $0 < c < \min \{
 |a|, b \}$ and that for $\vep$ satisfying $- \f{\varphi (-c)}{c} <
 \vep < \f{\varphi (c)}{c}$,   the infimum of $h (s, \vep)$ with
 respect to $s \in (a, b)$ is attained at the unique root,  $\zeta = \zeta (\vep)$, of the
 equation $\varphi^\prime (s) = \vep$ with respect to $s \in (-c, c)$.
 Moreover, $\zeta = \zeta (\vep) = \f{\vep}{ \si^2 } -
\f{\nu \vep^2}{2 \si^6} + O(\vep^3)$ for $- \f{\varphi (-c)}{c} <
 \vep < \f{\varphi (c)}{c}$.

\eeL

 \bpf  Let
$\varphi^\prime(s), \; \varphi^{\prime \prime} (s), \; \varphi^{(3)}
(s)$, and $\varphi^{(4)} (s)$ denote the first, second, third, and
fourth derivatives of $\varphi(s)$, respectively.  Note that
\[
\varphi (0) = 0, \qqu \varphi^\prime (0) = 0, \qqu \varphi^{\prime
\prime} (0) = \si^2, \qqu \varphi^{(3)} (0) = \nu.
\]
Clearly, $h(s, \vep)$ is convex with respect to $s \in (a, b)$. By
the assumption on $\vep$ and the convexity of $\varphi (s)$, we have
 \[
\lim_{ s \downarrow a } \varphi^\prime (s) \leq \lim_{ s \downarrow
a } \f{ \varphi
 (s)} {s} < \vep < \lim_{ s \uparrow b }  \f{ \varphi
 (s)} {s} \leq \lim_{ s \uparrow b } \varphi^\prime (s).
 \]
Since $\varphi (s)$ is convex for $s \in (a, b)$, it follows that
$\varphi^\prime (s)$ is increasing for $s \in (a, b)$.  Hence, the
equation $\varphi^\prime (s) = \vep$ with respect to $s \in (a, b)$
has a unique root,  $\zeta = \zeta (\vep)$, for $\vep$ such that
$\lim_{ s \downarrow a } \f{ \varphi
 (s)} {s} < \vep < \lim_{ s \uparrow b }  \f{ \varphi
 (s)} {s}$.  Moreover,
\[
\zeta =  \bec \zeta (\vep)   < 0 & \tx{for} \; \lim_{ s \downarrow a
} \f{ \varphi (s)} {s} < \vep < 0,\\
 \zeta (\vep)  = 0 & \tx{for} \; \vep = 0,\\
 \zeta (\vep)  > 0 & \tx{for} \; 0 < \vep < \lim_{ s \uparrow b } \f{ \varphi (s)} {s}
 \eec
\]
Since $h(s, \vep)$ is convex with respect to $s \in (a, b)$, it
follows that
 \[
\inf_{s \in (a, b)} h (s, \vep) = h (\zeta, \vep).
 \]
Note that $\varphi^{\prime \prime} (0) = \si^2 > 0$.  Since
$\varphi^{\prime \prime} (s)$ is a continuous function of $s \in (a,
b)$, it follows that there exists a positive number $c$ such that $0
< c < \min \{ |a|, \; b \}$ and that $\varphi^{\prime \prime} (s)
\geq \f{\si^2}{2}$ for all $s \in [ - c, c ]$.  By Taylor series
expansion formula,
\[
\varphi^\prime (s)  = \varphi^\prime (0) + \varphi^{\prime \prime}
(\eta s) s = \varphi^{\prime \prime} (\eta s) s \geq \f{\si^2}{2} s
\qu \tx{for $s \in (0, c]$,}
\]
where $\eta \in (0, 1)$ is a number dependent on $s$.
 Similarly, $\varphi^\prime (s) \leq \f{\si^2}{2} s$
for $s \in [-c, 0)$.   As a consequence of the convexity of $\varphi
(s)$, it must be true that $\varphi^\prime (- c) < \vep <
\varphi^\prime (c)$ for $\vep$ satisfying $- \f{\varphi (-c)}{c} <
 \vep < \f{\varphi (c)}{c}$.   Since $\varphi^\prime (s)$ is monotonically increasing with respect to
$s \in [- c, c ]$, it follows that $ - c < \zeta = \zeta (\vep) < c$
for $- \f{\varphi (-c)}{c} < \vep < \f{\varphi (c)}{c}$.   From now
on, we restrict $\vep$ to satisfy the
 constraint $- \f{\varphi (-c)}{c} <
 \vep < \f{\varphi (c)}{c}$.  Therefore,
\[
\vep = \varphi^\prime (\zeta) \geq \f{\si^2}{2} \zeta \qu \tx{for}
\qu 0 < \vep < \f{\varphi (c)}{c}
\]
and \[ \vep = \varphi^\prime (\zeta) \leq \f{\si^2}{2} \zeta \qu
\tx{for} \qu - \f{\varphi (-c)}{c} < \vep < 0.
\]
Hence, $\li | \f{ \zeta }{\vep} \ri | = \f{ \zeta }{\vep} \leq \f{ 2
}{ \si^2 }$ for $0 < \vep < \f{\varphi (c)}{c}$ and $- \f{\varphi
(-c)}{c} < \vep < 0$. This shows that \be \la{result88a}
 \zeta = \zeta (\vep) =
O (\vep) \ee for $- \f{\varphi (-c)}{c} < \vep < \f{\varphi
(c)}{c}$. By Taylor series expansion formula,
\[
\varphi^\prime (s) =  \si^2 s + \f{1}{2} \varphi^{(3)} (0) s^2 +
\f{1}{6} \varphi^{(4)} (\eta s) s^3 \qu \tx{for $s \in [- c, c]$,}
\]
where $\eta \in (0, 1)$.  Hence,
\[
\varphi^\prime (\zeta) =  \si^2 \zeta + \f{1}{2} \varphi^{(3)} (0)
\zeta^2 + \f{1}{6} \varphi^{ (4)} (\eta \zeta) \zeta^3
\]
for $- \f{\varphi (-c)}{c} < \vep < \f{\varphi (c)}{c}$,  where
$\eta \in (0, 1)$ is dependent on $\vep$.   Since $\varphi^{ (4) }
(s)$ is continuous with respect to $s \in [- c, c]$, there exists $K
> 0$ such that $| \varphi^{ (4) } (s) | \leq K$ for all $s \in [ - c, c
]$. Recall that  $h(s, \vep)$ is minimized at $s = \zeta = \zeta
(\vep)$ such that $\varphi^\prime (\zeta) = \vep$.   Hence, $\zeta$
satisfies the equation
\[
\si^2 \zeta + \f{1}{2} \varphi^{(3)} (0) \zeta^2 + \f{1}{6}
\varphi^{(4)} (\eta \zeta) \zeta^3 = \vep
\]
and thus \be \la{result88b}
 \zeta = \f{\vep}{ \si^2 } - \f{1}{2
\si^2} \varphi^{(3)} (0) \zeta^2 - \f{1}{6 \si^2} \varphi^{(4)}
(\eta \zeta) \zeta^3 \ee for $- \f{\varphi (-c)}{c} < \vep <
\f{\varphi (c)}{c}$.  Since $| \varphi^{(4)} (\eta \zeta) | \leq K$,
it follows from (\ref{result88a}) and (\ref{result88b}) that \bel
\zeta & = & \f{\vep}{ \si^2 } - \f{1}{2 \si^2} \varphi^{(3)} (0) [ O
(\vep) ]^2 - \f{1}{6 \si^2}
\varphi^{(4)} (\eta \zeta) [ O (\vep) ]^3 \nonumber\\
& = &  \f{\vep}{ \si^2 } + O (\vep^2) \la{result88c} \eel for $\vep$
satisfying $- \f{\varphi (-c)}{c} < \vep < \f{\varphi (c)}{c}$.
Again, since $| \varphi^{(4)} (\eta \zeta) | \leq K$, it follows
from (\ref{result88b}) and (\ref{result88c}) that  \bee \zeta & = &
\f{\vep}{ \si^2 } - \f{1}{2 \si^2} \varphi^{(3)} (0) \li [ \f{\vep}{
\si^2 } + O (\vep^2) \ri ]^2 - \f{1}{6 \si^2}
\varphi^{ (4)} (\eta \zeta) \li [ \f{\vep}{ \si^2 } + O (\vep^2) \ri ]^3\\
& = &  \f{\vep}{ \si^2 } - \f{ \nu \vep^2}{2 \si^6} + O(\vep^3) \eee
for $\vep$ satisfying $- \f{\varphi (-c)}{c} < \vep < \f{\varphi
(c)}{c}$. This completes the proof of the lemma.

\epf

\bsk

We are now in a position to prove the theorem.  The inequality
(\ref{firstinequality}) immediately follows from Theorem
\ref{special}. The inequality (\ref{secondinequality}) follows from
(\ref{firstinequality}) and the convexity of $\varphi(s)$. It
remains to investigate the asymptotic expression of the probability
bound $\li [ \exp  ( \varphi (\zeta) - \vep \zeta ) \ri ]^m$  and
the ratio $\f{\varphi(\zeta)}{\zeta}$.

Let $c$ be the number in the context of Lemma \ref{Order3388}. Using
Taylor series expansion formula, we have
\[ \varphi (s) = \varphi (0) + \varphi^\prime (0) s + \f{1}{2} \varphi^{\prime
\prime} (0) s^2 + \f{1}{6} \varphi^{(3)} (0) s^3 + O (s^4) =
\f{\si^2}{2} s^2 + \f{ \nu }{6} s^3 + O (s^4) \qu \tx{for $s \in [-
c, c]$. }
\]
Making use of this expression of $\varphi (s)$ and the expression of
$\zeta$ in Lemma \ref{Order3388}, we have \bee  \varphi (\zeta) & =
& \f{\si^2}{2} \zeta^2 + \f{ \nu }{6}
\zeta^3 + O (\zeta^4)\\
&  = & \f{\si^2}{2} \li [ \f{\vep}{ \si^2 } - \f{\nu \vep^2}{2
\si^6}
+ O(\vep^3)  \ri ]^2 + \f{ \nu }{6} \li [  \f{\vep}{ \si^2 } - \f{\nu \vep^2}{2 \si^6} + O(\vep^3) \ri ]^3 +  O( \vep^4 ) \\
& = & \f{\vep^2}{2 \si^2} - \f{\nu \vep^3}{2
\si^6} + \f{\nu \vep^3}{6 \si^6} + O (\vep^4)\\
& = &  \f{\vep^2}{2 \si^2} - \f{\nu \vep^3}{3 \si^6} + O (\vep^4),
\eee
\bee h (\zeta, \vep) & = & - \vep \zeta +
\varphi (\zeta)  \\
&  = &  - \f{\vep^2}{ \si^2 } + \f{\nu \vep^3}{2 \si^6} + O(\vep^4)
+ \f{\vep^2}{2 \si^2} - \f{\nu \vep^3}{3 \si^6} + O (\vep^4)\\
& = &  - \f{\vep^2}{2 \si^2} + \f{\nu \vep^3}{6 \si^6} + O (\vep^4),
\eee \bee \li [   \exp  ( \varphi (\zeta) - \vep \zeta ) \ri ]^m & =
& \li [ \exp
( h (\zeta, \vep)  ) \ri ]^m\\
& = & \li [  \exp \li (  - \f{\vep^2}{2 \si^2} + \f{\nu \vep^3}{6 \si^6} + O (\vep^4) \ri ) \ri ]^m\\
& = & \li [  \exp \li (  - \f{\vep^2}{2 \si^2} + \f{\nu \vep^3}{6
\si^6}  \ri )  \exp \li (  O
(\vep^4) \ri )\ri ]^m\\
& = & \li [  \exp \li (  - \f{\vep^2}{2 \si^2} + \f{\nu \vep^3}{6
\si^6}  \ri )  \ri ]^m \li [   \exp \li (  O
(\vep^4) \ri )\ri ]^m\\
& = & \li [  \exp \li (  - \f{\vep^2}{2 \si^2} + \f{\nu \vep^3}{6
\si^6}  \ri )  \ri ]^m \li [  1 +  O (\vep^4)  \ri ]^m\\
& = & \li [  1 +  O (\vep^4)  \ri ] \li [  \exp \li (  -
\f{\vep^2}{2 \si^2} + \f{\nu \vep^3}{6 \si^6}  \ri )  \ri ]^m\\
& = & \li [  1 +  O (\vep^3)  \ri ] \li [  \exp \li (  -
\f{\vep^2}{2 \si^2} \ri )  \ri ]^m,
 \eee
  and  \bee
\f{\varphi (\zeta)}{\zeta} & = & \f{ \f{\vep^2}{2 \si^2} - \f{\nu
\vep^3}{3 \si^6} + O (\vep^4) }{
\f{\vep}{ \si^2 } - \f{\nu \vep^2}{2 \si^6} + O(\vep^3) } \\
& = & \f{\vep}{2} \f{ \f{\vep}{ \si^2} - \f{2 \nu \vep^2}{3 \si^6} +
O (\vep^3) }{ \f{\vep}{ \si^2 } - \f{\nu \vep^2}{2 \si^6} +
O(\vep^3) }
\\
& = &  \f{\vep}{2} \f{ 1 - \f{2 \nu \vep}{3 \si^4} + O (\vep^2) }{1
- \f{\nu \vep}{2 \si^4} + O(\vep^2) }
\\
& = &  \f{\vep}{2}  \li [ 1 - \f{2 \nu \vep}{3 \si^4} + O (\vep^2)
\ri
]  \li [   1 + \f{\nu \vep}{2 \si^4} + O(\vep^2) \ri ] \\
& = &  \f{\vep}{2}  \li [ 1 - \f{\nu \vep}{6 \si^4} + O (\vep^2) \ri
]\\
& = &  \f{\vep}{2}  \li [ 1 +  O (\vep) \ri ] \eee for $\vep$
satisfying $- \f{\varphi (-c)}{c} < \vep < \f{\varphi (c)}{c}$.
This completes the proof of the theorem.

\section{Proof of Theorem \ref{Hoeffdingmu}} \la{Hoeffdingmuapp}

Define $h (\mu, s) = \ln ( \mu e^s + 1 - \mu )$ for $\mu \in [0, 1]$
and $s \in \bb{R}$. It is shown by Hoeffding in \cite{Hoeffding}
that
\[ \ln \bb{E}
\li [ \exp \li ( s \sum_{i = 1}^n X _i \ri ) \ri ] \leq \sum_{i =
1}^n h (\mu_i, s) \leq \mcal{V} (s, n)
\]
for all $s \in \bb{R}$ and $n \in \bb{N}$.  Let $\ga = m \se$.  Note
that $\mcal{V} (s, m) - \ga s = m \li [  h (\ovl{\mu}_m, s) - \se s
\ri ]$.  By differentiation, it can be readily shown that the
infimum of $h (\ovl{\mu}_m, s) - \se s$ with respect to $s \in
\bb{R}$ is attained at $\zeta$  and accordingly,
\[
\exp \li ( \mcal{V} (\zeta, m) - \ga \zeta \ri ) = \li [ \exp \li (
\se \ln \f{\ovl{\mu}_m}{\se} + (1 - \se) \ln \f{ 1 - \ovl{\mu}_m }{1
- \se} \ri ) \ri ]^m.
\]
Finally, invoking Theorem \ref{simple} yields the conclusion of the
theorem.

\section{Proof of Theorem \ref{Hoeffdingvar} } \la{Hoeffdingvarapp}

With the independence of the random variables,  it is shown by
Hoeffding in \cite{Hoeffding} that \[ \ln \bb{E} \li [ \exp \li ( s
\sum_{i = 1}^n X _i \ri ) \ri ] \leq \sum_{i = 1}^n \li [
\f{b^2}{b^2 + \si_i^2} \exp \li ( - \f{\si_i^2}{b} s \ri ) +
\f{\si_i^2}{b^2 + \si_i^2} e^{b s} \ri ] \leq \mcal{V} (s, n)
\]
for all $s \in \bb{R}$ and $n \in \bb{N}$. Let $\ga = m \vep$. Note
that $\mcal{V} (s, m) - \ga s = m \ln \li [ \f{b^2}{b^2 + \nu_m}
\exp \li (  - \f{\nu_m}{b} s \ri ) + \f{\nu_m}{b^2 + \nu_m} e^{b s}
\ri ] - m \vep s$. By differentiation, it can be readily shown that
the infimum of $\mcal{V} (s, m) - \ga s$ with respect to $s \in
\bb{R}$ is attained at $\zeta$ and accordingly,
\[
\exp \li ( \mcal{V} (\zeta, m) - \ga \zeta \ri ) = \li [ \li ( 1 +
\f{ b \vep }{ \nu_m } \ri )^{ - \f{ \nu_m + b \vep  }{ b^2 + \nu_m }
} \; \li (  1 - \f{\vep}{b} \ri )^{ - \f{ b^2 - b \vep }{ b^2 +
\nu_m } } \ri ]^m.
\]
Finally, invoking Theorem \ref{simple} yields the conclusion of the
theorem.

\section{Proof of Theorem \ref{Umiformnormal} } \la{Umiformnormalapp}

By the independence of the random variables,  \[ \ln \bb{E} \li [
\exp \li ( s \sum_{i = 1}^n X _i \ri ) \ri ] = \mcal{V} (s, n)
\]
for all $s \in \bb{R}$ and $n \in \bb{N}$.  Let $\ga = m \se$.  Note
that $\mcal{V} (s, m) - \ga s = m \li (  \ovl{\mu}_m s + \f{
\ovl{\nu}_m s^2}{2} - \se s  \ri )$.  By differentiation, it can be
readily shown that the infimum of $\ovl{\mu}_m s + \f{ \ovl{\nu}_m
s^2}{2} - \se s $ with respect to $s \in \bb{R}$ is attained at
$\zeta$ and accordingly,
\[
\exp \li ( \mcal{V} (\zeta, m) - \ga \zeta \ri ) = \li [ \exp \li (
- \f{ | \se - \ovl{\mu}_m |^2 }{ 2 \ovl{\nu}_m} \ri )  \ri ]^m.
\]
Finally, invoking Theorem \ref{simple} yields the conclusion of the
theorem.

\section{Proof of Theorem \ref{UniformPoisson} } \la{UniformPoissonapp}

By the independence of the random variables,  \[ \ln \bb{E} \li [
\exp \li ( s \sum_{i = 1}^n X _i \ri ) \ri ] = \mcal{V} (s, n)
\]
for all $s \in \bb{R}$ and $n \in \bb{N}$.  Let $\ga = m \se$.  Note
that $\mcal{V} (s, m) - \ga s = n \ovl{\lm}_m (e^s - 1)  - m \se s$.
By differentiation, it can be readily shown that the infimum of
$\mcal{V} (s, m) - \ga s$ with respect to $s \in \bb{R}$ is attained
at $\zeta$ and accordingly,
\[
\exp \li ( \mcal{V} (\zeta, m) - \ga \zeta \ri ) = \li [ \exp \li (
\se - \ovl{\lm}_m + \se \ln \f{\ovl{\lm}_m}{\se} \ri ) \ri ]^m.
\]
Finally, invoking Theorem \ref{simple} yields the conclusion of the
theorem.

\end{document}